\documentclass{IEEEtran4PSCC}
\ifCLASSINFOpdf
   \usepackage[pdftex]{graphicx}
\else
   \usepackage[dvips]{graphicx}
\fi
\usepackage{bm}
\usepackage{multirow}
\usepackage{multicol}
\usepackage{cite}

\usepackage[cmex10]{amsmath}
\usepackage{mathtools}

\DeclarePairedDelimiter\floor{\lfloor}{\rfloor}

\newtheorem{thm}{Theorem}[section]
\newtheorem{lemma}{Lemma}[section]

\newtheorem{coro}{Corollary}[section]

\DeclareMathOperator{\Tr}{Tr}

\usepackage{amsmath,amssymb}

\usepackage{color}
\usepackage{amsfonts}
\DeclareMathAlphabet{\mymathbb}{U}{BOONDOX-ds}{m}{n} 

\DeclareMathOperator{\re}{Re}
\DeclareMathOperator{\im}{Im}
\DeclareMathOperator{\rank}{rank}

\makeatletter
\let\old@ps@headings\ps@headings
\let\old@ps@IEEEtitlepagestyle\ps@IEEEtitlepagestyle
\def\psccfooter#1{%
    \def\ps@headings{%
        \old@ps@headings%
        \def\@oddfoot{\strut\hfill#1\hfill\strut}%
        \def\@evenfoot{\strut\hfill#1\hfill\strut}%
    }%
    \def\ps@IEEEtitlepagestyle{%
        \old@ps@IEEEtitlepagestyle%
        \def\@oddfoot{\strut\hfill#1\hfill\strut}%
        \def\@evenfoot{\strut\hfill#1\hfill\strut}%
    }%
    \ps@headings%
}
\makeatother

\begin{document}
%
\title{Balancibility: Existence and Uniqueness of Power Flow Solutions under Voltage Balance Requirements}

\author{
\IEEEauthorblockN{Bowen Li, Bai Cui, and Feng Qiu}
\IEEEauthorblockA{Energy Systems Division\\Argonne National Laboratory\\
Lemont, IL, USA
\\
\{bowen.li, bcui, fqiu\}@anl.gov}
\and
\IEEEauthorblockN{Daniel K. Molzahn}
\IEEEauthorblockA{School of Electrical and Computer Engineering\\
 Georgia Institute of Technology,\\
 Atlanta, GA, USA\\molzahn@gatech.edu}
}


\maketitle

\begin{abstract}
In distribution systems, power injection variability due to growing penetrations of distributed energy resources (DERs) and dispatchable loads can lead to power quality issues such as severe voltage unbalance. To ensure safe operation of phase-balance-sensitive components such as transformers and induction motor loads, the amount of voltage unbalance must be maintained within specified limits for a range of uncertain loading conditions. This paper builds on existing ``solvability conditions'' that characterize operating regions for which the power flow equations are guaranteed to have a unique high-voltage solution. We extend these existing solvability conditions to be applicable to distribution systems and augment them with a ``balancibility'' condition which quantifies an operating region within which a unique, adequately balanced power flow solution exists. To build this condition, we consider different unbalance definitions and derive closed-form representations through reformulations or safe approximations. Using case studies, we evaluate these closed-form representations and compare the balancibility conditions associated with different unbalance definitions.
\end{abstract}




\section{Introduction}
Increasing penetrations of distributed energy resources (DERs) and dispatchable loads can result in greater variability and stochasticity of the power injections in distribution systems. Extreme variations in power injections can also lead to power quality issues such as significant voltage unbalances in which the voltage magnitudes and angles have large offsets between the three phases. Unbalanced voltages can greatly impact essential power system devices such as three-phase induction motors and transformers~\cite{bal_motor1,bal_motor2}. Specifically, for induction motors, even small amounts of voltage unbalance can cause severe temperature rise, efficiency loss, and decreased life expectancy, which leads to serious consequences from premature motor failures, costly shutdowns, and lost production \cite{doe_bal1}. Voltage unbalance can cause annual losses to U.S. industries of up to $\$28$ billion~\cite{doe_bal2}. Hence, it is critical to provide secure criteria for power system operations subject to uncertainties (e.g., renewable generation or load consumption) such that the resulting steady-state operating points have voltage balance guarantees.

Organizations such as International Electrotechnical Commission (IEC), National Electrical Manufacturers Association (NEMA), and IEEE have each developed definitions to quantify the amount of voltage unbalance. For example, IEC~\cite{iec} and IEEE~\cite{IEEEstd2} have definitions that are based on the ratio between negative/zero-sequence voltage and positive-sequence voltage calculated from the symmetrical component transformation. Othe standards from NEMA~\cite{nema} and IEEE~\cite{ieeestd, ieeestd3} define voltage unbalance using line-to-line and line-to-ground voltage magnitudes, respectively. Previous works \cite{zero_seq,ktij_bal} summarize and compare these definitions. 

This paper characterizes regions of power injections for which the power flow equations admit a unique high-voltage power flow solution that satisfies specified phase unbalance requirements according to these definitions. 
The approach in this paper builds on existing power flow solvability conditions, which have been extensively studied for both transmission systems~\cite{lee2018convex, dji_sol, cui2019solvability} and distribution systems~\cite{bolo_sol, wang_sol, wang_sol2, andrey_loadflow}. This paper extends existing solvability conditions to consider voltage balance requirements, resulting in the proposed ``balancibility'' condition. This balancibility condition quantifies a region of power injections for which a unique and balanced high-voltage power flow solution is guaranteed to exist. 

To the best of our knowledge, this is the first paper to incorporate voltage balance requirements into power flow solvability conditions. As specific contributions, we consider different unbalance definitions and develop various approaches for deriving closed-form reformulations or safe approximations that quantify the voltage unbalance level. We say a set is a ``safe approximation'' of a robust set if it contains the entire robust set. We use a general model to describe the sets that contain the power flow solutions under uncertain power injections and provide supporting theoretical guarantees on the quality of the approaches. We demonstrate the proposed balancibility condition using the solvability condition in \cite{cui2019solvability}. We then numerically illustrate the quality of the balancibility conditions associated with different unbalance definitions. 

The proposed balancibility condition is expected to be a key enabling tool for many applications due to its ability to greatly simplify various problem formulations, particularly those which require voltage balance guarantees under uncertain power injections (e.g., robust AC optimal power flow). This condition is also useful for identifying the worst-case uncertainty realizations with respect to voltage balance limits. 

The remainder of the paper is organized as follows. Section~\ref{sec:not} introduces notation. Section~\ref{sec:model} describes the distribution network model and the solvability condition from~\cite{cui2019solvability}. Section~\ref{sec:bal} derives the closed-form approximations or reformulations for various unbalance definitions and proposes our balancibility condition. Section~\ref{sec:case} presents case studies. Section~\ref{sec:conclu} summarizes the paper and discusses future directions.

\section{Notation}\label{sec:not}
Boldface letters indicate complex variables and roman font is used for real variables. $\bm{j}=\sqrt{-1}$. Transposition and Hermitian transposition are denoted as $(\cdot)^{\top}$ and $(\cdot)^{\text{H}}$, respectively. $\mathbb{I}_n\in\mathbb{R}^{n\times n}$ represents the identity matrix. $\mymathbb{0}_n$ denotes an $n\times n$ zero matrix. $\mathcal{H}^n$ denotes the set of $n\times n$ Hermitian matrices. For $y\in\mathbb{R}$, $\floor{y}$ returns the greatest integer less than or equal to $y$. For $x\in\mathbb{R}^n$ or $\mathbb{C}^n$, $\bar{x}$ denotes its component-wise conjugate. $\re(\bm{x})$ and $\im(\bm{x})$ denote its component-wise real and imaginary parts, respectively. $x_i$ denotes the $i$-th entry in the vector and $x_{i,j}$ ($i\leq j$) denotes the vector from $i$-th entry to $j$-th entry. $\|x\|$ denotes the $\ell_2$-norm and $\|x\|_p$ denotes the $\ell_p$-norm. $|\bm{x}|$ returns the magnitude if $\bm{x} \in \mathbb{C}$. All angle values are reported in degrees. $\mathcal{A}\times \mathcal{B}$ denotes the Cartesian product of sets $\mathcal{A}$ and $\mathcal{B}$. $\mathcal{A}^n$ represents the Cartesian product with set $\mathcal{A}$ for $n$ times. For matrix $X\in\mathbb{C}^{n\times n}$ or $\mathbb{R}^{n\times n}$, $X_{ij}$ represents the entry at $i$-th row and $j$-th column. $X_i$ represents the vector of the $i$-th row. $X_i^{j,k}$ ($j\leq k$) represents the vector from \mbox{$j$-th} to \mbox{$k$-th} elements in the $i$-th row of $X$. $N(X)$ represents the nullspace of $X$. $\rank(X)$ and $\Tr(X)$ return the rank and trace, respectively, of $X$. $X\succeq 0$ indicates positive semidefiniteness of $X$. The function $\mathbf{\lambda}_{\min}(X)$ returns the smallest eigenvalue of $X$. The function $\text{blkdiag}(\cdot)$ returns a block diagonal matrix with its input matrices and $\text{diag}(\cdot)$ returns a diagonal matrix. $\mathbf{0}$ ($\mathbf{1}$) represents all-zero (all-one) vector or matrix with appropriate size. For a set $S$, its closure and boundary are denoted by $\bar{S}$ and $\partial S$, respectively. $D(x,r)$ represents an open disk with center $x$ and radius $r$. For brevity, we denote $D(0,r)$ by $D(r)$.


\section{Network model and solvability Condition}\label{sec:model}
In this paper, we use a distribution network model similar to \cite{andrey_loadflow,cui2019solvability} and assume a generic network topology (i.e., radial or meshed) with a single slack bus and multiple-phase wye-connected PQ buses. We choose the slack bus to be at node $0$ and define $\mathcal{N}_L$ as the set of PQ buses. Denote the voltage at the slack bus as $\mathbf{V}_G=(\mathbf{V}_{G,a},\mathbf{V}_{G,b},\mathbf{V}_{G,c})^{\top}$ for each phase $a,b,c$. Similarly, for all $i\in\mathcal{N}_L$, we define its wye-connected power consumption to be $\mathbf{S}_L^i=(\mathbf{S}^i_{L,a},\mathbf{S}^i_{L,b},\mathbf{S}^i_{L,a})^{\top}$ and its voltage to be $\mathbf{V}_L^i=(\mathbf{V}^i_{L,a},\mathbf{V}^i_{L,b},\mathbf{V}^i_{L,a})^{\top}$. Based on the admittance matrix $\mathbf{Y}$, we have 
\begin{align}
    \begin{bmatrix}\mathbf{I}_G\\-\mathbf{I}_L\end{bmatrix}=\begin{bmatrix}\mathbf{Y_{GG}} & \mathbf{Y_{GL}}\\\mathbf{Y_{LG}} & \mathbf{Y_{LL}}\end{bmatrix}\begin{bmatrix}\mathbf{V}_G\\-\mathbf{V}_L\end{bmatrix}
\end{align}
where $\mathbf{I}_G$ is the current injected at the slack bus and $\mathbf{I}_L$ is the current withdrawn at PQ buses.\footnote{For two-phase or single-phase nodes, $\mathbf{V}_L$, $\mathbf{S}_L$, and $\mathbf{Y}$ only collect quantities for the existing phases.} Based on \cite{andrey_loadflow,cui2019solvability}, we have 
\begin{align}
    \mathbf{v}_L=\mathbf{1}-\hat{\mathbf{Z}}\text{diag}^{-1}(\bar{\mathbf{v}}_L)\bar{\mathbf{S}}_L
\end{align}
where\footnote{The invertibility of $\mathbf{Y}_{LL}$ is proved in \cite{wang_sol}. As in \cite{andrey_loadflow,wang_sol,cui2019solvability}, we assume that $\mathbf{E}$ and $\mathbf{v}_L$ do not contain zero elements, which is the case for practical power systems.}
\begin{subequations}
\begin{align}
\mathbf{E} &= -\mathbf{Y}_{LL}^{-1}\mathbf{Y}_{LG}\mathbf{V}_{G},\\ \mathbf{v}_L &= \text{diag}^{-1}(\mathbf{E})\,\mathbf{V}_{L}, \\
\hat{\mathbf{Z}} &= \text{diag}^{-1}(\mathbf{E})\,\mathbf{Y}_{LL}^{-1}\,\text{diag}^{-1}(\bar{\mathbf{E}}).
\end{align}
\end{subequations}

\begin{figure}
\centering
\vspace{0cm}
\includegraphics[width=3.5in]{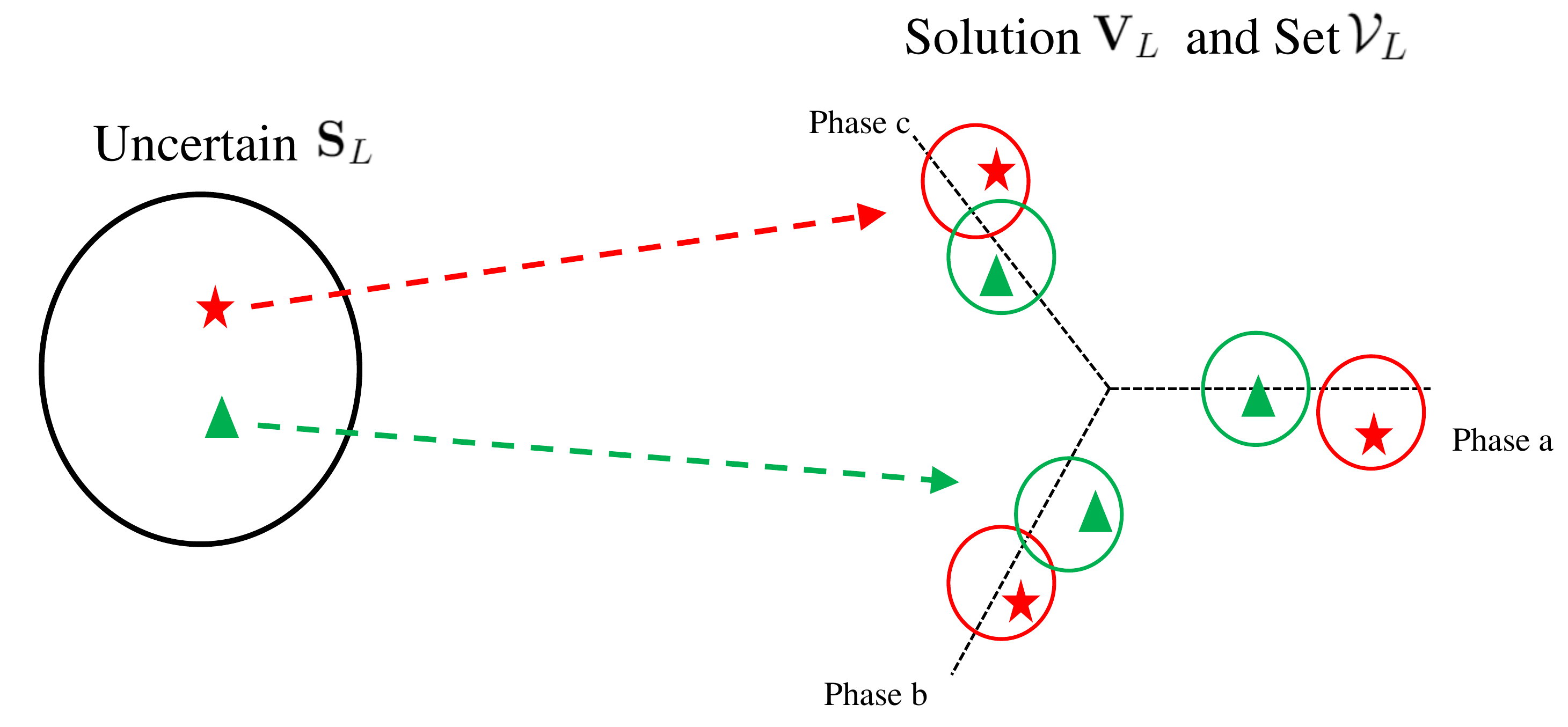}
\caption{Relationship among $\mathbf{S}_L$, $\mathbf{V}_L$, and $\mathcal{V}_L$}
\vspace{0cm}
\label{fig: fig6}
\end{figure}

We use the solvability condition from \cite{cui2019solvability} to analyze the secure region of $\mathbf{S}_L$ for which there exists a unique $\mathbf{V}_L$ within a set $\mathcal{V}_L(\mathbf{S}_L)$ (i.e., parameterized on $\mathbf{S}_L$). As illustrated in Fig.~\ref{fig: fig6}, if $\mathbf{S}_L$ changes in its uncertainty set, $\mathcal{V}_L$ also changes to follow the power flow solution $\mathbf{V}_L$ if the solvability condition is satisfied. Since the variation of $\mathcal{V}_L$ can be easily represented as an explicit function of $\mathbf{S}_L$, the solvability condition summarized in this section provides an efficient way to confidently locate the power flow solution under uncertainty. A detailed description of this condition is available in~\cite{cui2019solvability}. 

Define a nominal power flow solution $(\mathbf{v}^0_L,\mathbf{S}^0_L)$,
\begin{align}
\mathbf{v}^0_L=\mathbf{1}-\hat{\mathbf{Z}}\,\text{diag}^{-1}(\bar{\mathbf{v}}^0_L)\,\bar{\mathbf{S}}^0_L,    
\end{align}
and $\bm{\sigma}_L=\mathbf{S}_L-\mathbf{S}^0_L$. If no nominal solution is provided, a trivial selection is $\mathbf{v}^0_L=\mathbf{1}$ when $\mathbf{S}^0_L=0$. We also define the following quantities used in the solvability condition and $\mathcal{V}_L$:
\begin{align*}
 &\tilde{\mathbf{Z}}=\text{diag}^{-1}(\mathbf{v}_L^0)\,\hat{\mathbf{Z}}\,\text{diag}^{-1}(\bar{\mathbf{v}}_L^0),\ \mathbf{u}_L=\text{diag}^{-1}(\mathbf{v}^0_L)\,\mathbf{v}_{L}.   
\end{align*}

For $i\in\mathcal{N}_L$ and phase $p\in\{a,b,c\}$, denote $\tilde{\mathbf{Z}}_i^p$ as the corresponding row of $\tilde{\mathbf{Z}}$, then
\begin{align*}
    &\eta_{i,p}(\bm{\sigma}_L)=(\tilde{\mathbf{Z}}_i^p)^{\top}\bar{\bm{\sigma}}_L,\quad \xi_{i,p}(\mathbf{S}_L)=\|(\tilde{\mathbf{Z}}_i^p)^{\top}\text{diag}(\bar{\mathbf{S}}_L)\|_1,\\
    \nonumber & \gamma_{i,p}(\bm{\sigma}_L,\mathbf{S}_L)=2(\xi_{i,p}(\mathbf{S}_L)+\re(\eta_{i,p}(\bm{\sigma}_L)))\\
    &\hspace{1.5in}-|\eta_{i,p}(\bm{\sigma}_L)|^2-\xi_{i,p}(\mathbf{S}_L)^2,
\end{align*}
where $\eta_{i,p}(\bm{\sigma}_L)$ and $\xi_{i,p}(\mathbf{S}_L)$ represent aggregated system stress measures on each node and phase resulting from incremental and total loads. These measures also appear in other existing solvability literature \cite{dji_sol0,wang_sol,dji_sol}. We also define $\gamma_{i,p}(\bm{\sigma}_L,\mathbf{S}_L)$ fusing these two stresses. Accordingly, we define 
\begin{align*}
&\eta(\bm{\sigma}_L)=\max_{i\in\mathcal{N}_L,\ p\in\{a,b,c\}}|\eta_{i,p}(\bm{\sigma}_L)|, \\
&\xi(\mathbf{S}_L)=\max_{i\in\mathcal{N}_L,\ p\in\{a,b,c\}}\xi_{i,p}(\mathbf{S}_L),\\
&\gamma(\bm{\sigma}_L,\mathbf{S}_L)=\max_{i\in\mathcal{N}_L,\ p\in\{a,b,c\}}\gamma_{i,p}(\bm{\sigma}_L,\mathbf{S}_L),\\
&\Delta=(1-\gamma(\bm{\sigma}_L,\mathbf{S}_L))^2-4\xi(\mathbf{S}_L)^2\eta(\bm{\sigma}_L)^2.
\end{align*}

Next, we construct the following framework which geometrically quantifies a disk for $\mathbf{u}_L$ with parameter $r\geq0$. For $i\in\mathcal{N}_L$ and $p\in\{a,b,c\}$, we have
\begin{subequations}\label{eq: invset}
\begin{align}
&\xi_{i,p}(\mathbf{S}_L)>0:\quad |1-\eta_{i,p}(\bm{\sigma}_L)-\mathbf{u}_{L,p}^i|\leq r\,\xi_{i,p}(\mathbf{S}_L),\\
&\xi_{i,p}(\mathbf{S}_L)=0:\quad \mathbf{u}^i_{L,p}=1-\eta_{i,p}(\bm{\sigma}_L).
\end{align}
\end{subequations}

Hence, when $\xi_{i,p}(\mathbf{S}_L)=0$, $\mathbf{u}_{L,p}^i$ degenerates into a single point. The relationship between $r$ and these quantities is presented in the next theorem.
\begin{thm}\label{thm:solv}(Theorem B.3 in \cite{cui2019solvability}) Given nominal solution $(\mathbf{v}^0_L,\mathbf{S}^0_L)$, if the following condition is satisfied
\begin{subequations}\label{eq:solv}
\begin{align}
   &\gamma(\bm{\sigma}_L,\mathbf{S}_L)+2 \xi(\mathbf{S}_L)\eta(\bm{\sigma}_L)<1,\\
   &\xi(\mathbf{S}_L)-\eta(\bm{\sigma}_L)\leq 1.
\end{align}
\end{subequations}

Then, there exists a unique solution $\mathbf{u}_L$ in \eqref{eq: invset} with
\begin{equation}
\begin{cases}
r=\sqrt{\frac{1-\gamma(\bm{\sigma}_L,\mathbf{S}_L)-\sqrt{\Delta}}{2\xi(\mathbf{S}_L)^2}},\quad\mbox{if }\xi(\mathbf{S}_L)>0,\\
r=0,\quad\mbox{if }\xi(\mathbf{S}_L)=0.
\end{cases}    
\end{equation}
\end{thm}

From Theorem~\ref{thm:solv}, we know where $\mathbf{u}_L$ is located under $\mathbf{S}_L$, which can then be used to obtain $\mathcal{V}_L$. For any $i\in\mathcal{N}_L$ and $p=\{a,b,c\}$, we have $\mathbf{V}_{L,p}^i=\mathbf{E}_{p}^i\mathbf{v}_{L,p}^{0,i}\mathbf{u}_{L,p}^i$ where $\mathbf{E}_{p}^i$ and $\mathbf{v}_{L,p}^{0,i}$ are corresponding terms in $\mathbf{E}$ and $\mathbf{v}_L^0$. Then, the set $\mathcal{V}_{L,p}^i$ that contains $\mathbf{V}_{L,p}^i$ can be represented as
\begin{equation}
|(1-\eta_{i,p}(\bm{\sigma}_L))\mathbf{E}_{p}^i\mathbf{v}_{L,p}^{0,i}-\mathbf{V}_{L,p}^i|\leq r|\mathbf{E}_{p}^i\mathbf{v}_{L,p}^{0,i}|\xi_{i,p}(\mathbf{S}_L). \label{eq: invset2}  
\end{equation}

With $\mathcal{V}_L$, we then have the foundations to analyze voltage unbalance levels under $(\mathbf{V}_L,\mathbf{S}_L)$ under uncertainty and derive balancibility conditions that limit this level.   

\section{Safe Approximation on Balancibility}\label{sec:bal}

This section derives our proposed the balancibility conditions which guarantee satisfaction of the voltage balance requirements for all the realizations in $\mathcal{V}_L$. We use safe approximations and different reformulation techniques to develop closed-form representations of these conditions. We choose safe approximations since relaxations may underestimate the true voltage unbalance level and give insecure results. We assume that there are critical nodes ($i^*\in\mathcal{N}_L$) that are sensitive to amounts of voltage unbalance outside of specified limits. Then, for each critical node, we rewrite $\mathcal{V}_L$ as a general set $\mathbf{U}_{in}$ (with $\mathbf{U}_{in}^p$ denoting a particular phase $p=\{a,b,c\}$) as
\begin{equation}\label{eq: ginv}
    |\mathbf{V}_a-\mathbf{C}_a|\leq r_a,\ |\mathbf{V}_b-\mathbf{C}_b|\leq r_b,\ |\mathbf{V}_c-\mathbf{C}_a|\leq r_c
\end{equation}
where subscripts $a,b,c$ denote the phases. Geometrically, $\mathbf{U}_{in}^p$ is a disk $\bar{D}(\mathbf{C}_p,r_p)$ with center $\mathbf{C}_p$ and radius $r_p$. To help some of the derivations, we also represent these sets in real coordinates, denoted $\mathcal{V}_{in}$, by separating the real and imaginary parts in \eqref{eq: ginv}
\begin{subequations}\label{eq: ginv2}
\begin{align}
    &(V_a^r-C_a^r)^2+(V_a^i-C_a^i)^2\leq r_a^2,\label{eq: 1a}\\
    &(V_b^r-C_b^r)^2+(V_b^i-C_b^i)^2\leq r_b^2,\\
    &(V_c^r-C_c^r)^2+(V_c^i-C_c^i)^2\leq r_c^2,
\end{align}
\end{subequations}
where superscripts $r,i$ denotes the real and imaginary parts. Note that $\mathbf{U}_{in}$/$\mathcal{V}_{in}$ can be seen as a general output from any solvability conditions in complex domain \cite{andrey_loadflow, wang_sol, wang_sol2} or involving voltage magnitudes \cite{dji_sol}. Hence, the applicability is not restricted to any network assumptions (e.g., radial or wye-connected loads) or particular solvability condition \cite{cui2019solvability}. 

Linking back to solvability condition \eqref{eq: invset2} (taking phase $a$ as an example), $V_a^r$ and $V_a^i$ represent the real and imaginary parts of the complex voltage $\mathbf{V}_a=\mathbf{V}_{L,a}^{i^*}$. $C_a^r$ and $C_a^i$ represent the real and imaginary parts of $\mathbf{C}_a=(1-\eta_{i^*,a}(\bm{\sigma}_L))\mathbf{E}_{a}^{i^*}\mathbf{v}_{L,a}^{0,i^*}$ and $r_a=r|\mathbf{E}_{a}^{i^*}\mathbf{v}_{L,a}^{0,i^*}|\xi_{i^*,a}(\mathbf{S}_L)$. If $\xi_{i^*,a}(\mathbf{S}_L)=0$, \eqref{eq: invset} is degenerate and $V_a^r$ and $V_a^i$ can be treated as constants while analyzing the voltage unbalance and hence do not affect the results. 

To make \eqref{eq: ginv2} concise, we define vectors $V_a=(V_a^r, V_a^i)^{\top}\in\mathbb{R}^2$, $C_a=(C_a^r, C_a^i)^{\top}\in\mathbb{R}^2$, and set $\mathcal{V}^a_{in}\subset \mathbb{R}^2:=\{V_a \mbox{ that satisfies } \eqref{eq: 1a}\}$, with the same notations applied to phase $b$ and $c$. We also define $V_{abc}=(V_a^\top,V_b^\top,V_c^\top)^{\top}\in\mathbb{R}^6$, $\mathbf{V}_{abc}=(\mathbf{V}_a,\mathbf{V}_b,\mathbf{V}_c)^{\top}\in\mathbb{C}^3$, $r_{abc}=(r_a,r_b,r_c)^{\top}\in\mathbb{R}_+^3$, and $C_{abc}=(C_a^\top,C_b^\top,C_c^\top)^{\top}\in\mathbb{R}^6$.  Next, we derive the reformulations or safe approximations of the voltage balance requirement for different unbalance definitions. 
\subsection{Phase Voltage Unbalance Rate (PVUR) Definition}
In \cite{ieeestd}, the following definition of phase voltage unbalance rate ($PVUR$) is provided using the line-to-ground voltage magnitudes $|\mathbf{V}_a|$, $|\mathbf{V}_b|$, and $|\mathbf{V}_c|$:
\begin{equation}
    PVUR=\Delta^{\max}_V/V_{\mathrm{avg}}, \label{eq: PVUR}
\end{equation}
where $V_{\mathrm{avg}}=\frac{|\mathbf{V}_a|+|\mathbf{V}_b|+|\mathbf{V}_c|}{3}$ and 
$$\Delta^{\max}_V=\max\{||\mathbf{V}_a|-V_{\mathrm{avg}}|,||\mathbf{V}_b|-V_{\mathrm{avg}}|,||\mathbf{V}_c|-V_{\mathrm{avg}}|\}.$$

To ensure the power flow solutions are balanced, we require that the voltage profile satisfies \eqref{eq: PVUR} with a predefined tolerance of $\epsilon\in(0,1)$ with $PVUR\leq \epsilon$. This requirement is equivalent to the following linear constraints \cite{linepsccbal}:
\begin{align}
   \begin{bmatrix}
        \epsilon+2 & \epsilon-1 & \epsilon-1 \\ 
        \epsilon-1 & \epsilon+2 & \epsilon-1 \\     
        \epsilon-1 & \epsilon-1 & \epsilon+2 \\ 
        \epsilon-2 & \epsilon+1 & \epsilon+1 \\ 
        \epsilon+1 & \epsilon-2 & \epsilon+1 \\ 
        \epsilon+1 & \epsilon+1 & \epsilon-2 \\ 
     \end{bmatrix} 
     \begin{bmatrix}
        |\mathbf{V}_a| \\ 
        |\mathbf{V}_b| \\        
        |\mathbf{V}_c| \\ 
     \end{bmatrix} \geq \mathbf{0}.\label{eq: PVUR_lin}
\end{align}

Next, we require that all the solutions $V_{abc}\in \mathcal{V}_{in}$ satisfy \eqref{eq: PVUR_lin}. Without loss of generality, we only use the first linear constraint in \eqref{eq: PVUR_lin} as an example and the problem becomes
\begin{equation}
   \min_{V_{abc}\in\mathcal{V}_{in}} \big\{ (\epsilon+2)|\mathbf{V}_a| + (\epsilon-1)|\mathbf{V}_b|+ (\epsilon-1)|\mathbf{V}_c| \big\} \geq 0. \label{eq: pvurrob}
\end{equation}
Since $\mathcal{V}_{in}$ is separable in each phase, \eqref{eq: pvurrob} is equivalent to
\begin{multline}
    \min_{V_{a}\in\mathcal{V}^a_{in}} (\epsilon+2)|\mathbf{V}_a| + \min_{V_{b}\in\mathcal{V}^b_{in}} (\epsilon-1)|\mathbf{V}_b| \\
    + \min_{V_{c}\in\mathcal{V}^c_{in}} (\epsilon-1)|\mathbf{V}_c| \geq 0. \label{eq: pvurrob2}   
\end{multline}
Each subproblem in \eqref{eq: pvurrob2} can be easily solved since $\mathcal{V}^a_{in}$, $\mathcal{V}^b_{in}$, and $\mathcal{V}^c_{in}$ are closed disks. Since $\epsilon+2>0$ and $\epsilon-1<0$ ($\epsilon\in(0,1)$), we have the following reformulation of \eqref{eq: pvurrob2}:
\begin{align}
    (\epsilon+2)\max\{\|C_a\|-r_a,0\} + &(\epsilon-1)(r_b+\|C_b\|) \\
    + &(\epsilon-1)(r_c+\|C_c\|) \geq 0.
\end{align}
We use $\max\{\|C_a\|-r_a,0\}$ in case $|\mathbf{V}_a|=0$ when $\mathcal{V}^a_{in}$ contains the origin. Now, we can derive the voltage balance requirements using each of the linear constraints in \eqref{eq: PVUR_lin}. This approach also applies to other $PVUR$ definitions as in \cite{zero_seq,ieeestd3} where
$$\Delta^{\max}_V=\max\{|\mathbf{V}_a|,|\mathbf{V}_b|,|\mathbf{V}_c|\}-\min\{|\mathbf{V}_a|,|\mathbf{V}_b|,|\mathbf{V}_c|\}.$$

\subsection{Line Voltage Unbalance Rate (LVUR) Definition}
In \cite{nema}, an unbalance definition called the line voltage unbalance rate ($LVUR$) is provided using line-to-line voltages $|\mathbf{V}_{ab}|\!=\!|\mathbf{V}_{a}-\mathbf{V}_{b}|$, $|\mathbf{V}_{bc}|\!=\!|\mathbf{V}_{b}-\mathbf{V}_{c}|$, and $|\mathbf{V}_{ca}|\!=\!|\mathbf{V}_{c}-\mathbf{V}_{a}|$:
\begin{equation}
    LVUR=\Delta^{\max}_{V_L}/V_{\mathrm{avg},L}, \label{eq: LVUR}
\end{equation}
where $V_{\mathrm{avg},L}=\frac{|\mathbf{V}_{ab}|+|\mathbf{V}_{bc}|+|\mathbf{V}_{ca}|}{3}$ and 
$$\Delta^{max}_{V_L}=\max\{||\mathbf{V}_{ab}|-V_\mathrm{avg}|,||\mathbf{V}_{bc}|-V_\mathrm{avg}|,||\mathbf{V}_{ca}|-V_\mathrm{avg}|\}.$$

Similar to $PVUR$, with voltage balance requirement $LVUR\leq\epsilon$, we have
\begin{align}
   \begin{bmatrix}
        \epsilon+2 & \epsilon-1 & \epsilon-1 \\ 
        \epsilon-1 & \epsilon+2 & \epsilon-1 \\     
        \epsilon-1 & \epsilon-1 & \epsilon+2 \\ 
        \epsilon-2 & \epsilon+1 & \epsilon+1 \\ 
        \epsilon+1 & \epsilon-2 & \epsilon+1 \\ 
        \epsilon+1 & \epsilon+1 & \epsilon-2 \\ 
     \end{bmatrix} 
     \begin{bmatrix}
        |\mathbf{V}_{ab}| \\ 
        |\mathbf{V}_{bc}| \\        
        |\mathbf{V}_{ca}| \\ 
     \end{bmatrix} \geq \mathbf{0}.\label{eq: LVUR_lin}
\end{align}

We require that all $V_{abc}\in\mathcal{V}_{in}$ satisfy \eqref{eq: LVUR_lin}. Here, we use the first constraint in \eqref{eq: LVUR_lin} as an example:
\begin{equation}
   \left\{\min_{V_{abc}\in\mathcal{V}_{in}} (\epsilon+2)|\mathbf{V}_{ab}| + (\epsilon-1)|\mathbf{V}_{bc}|+ (\epsilon-1)|\mathbf{V}_{ca}|\right\}\geq 0. \label{eq: lvurrob}
\end{equation}

There are several approaches for safely approximating \eqref{eq: lvurrob}. The first approach bounds $|\mathbf{V}_{ab}|$, $|\mathbf{V}_{bc}|$, and $|\mathbf{V}_{ca}|$ as in \cite{LVURbound}. Taking $|\mathbf{V}_{ab}|$ as an example, we have
$$|\mathbf{V}_{ab}|=\|C_a-C_b+r_au_a+r_bu_b\|$$
where $u_a$ and $u_b$ are any unit vectors in $\mathbb{R}^2$. Then, we see that
\begin{subequations}
\begin{align}
&\max\{\|C_a-C_b\|-r_a-r_b,0\} \leq |\mathbf{V}_{ab}|,\\
& \|C_a-C_b\|+r_a+r_b\geq |\mathbf{V}_{ab}|.   
\end{align}
\end{subequations}

The voltages $|\mathbf{V}_{bc}|$ and $|\mathbf{V}_{ca}|$ are bounded analogously. Using a similar idea as in \eqref{eq: pvurrob2}, we safely approximate \eqref{eq: lvurrob} as
\begin{align}
    &(\epsilon+2)\max\{\|C_a-C_b\|-r_a-r_b,0\}\nonumber\\
    & \quad + (\epsilon-1)(\|C_b-C_c\|+r_b+r_c) \nonumber\\
    &\quad +(\epsilon-1)(\|C_a-C_b\|+r_b+r_c) \geq 0.\label{eq: llvbound}
\end{align}

The second approach for approximating \eqref{eq: lvurrob} uses the following relationship:
\begin{align*}
&\big||\mathbf{V}_a|-|\mathbf{V}_b|\big| \leq|\mathbf{V}_{ab}|\leq |\mathbf{V}_a|+|\mathbf{V}_b|,\\
& \big||\mathbf{V}_b|-|\mathbf{V}_c|\big| \leq|\mathbf{V}_{bc}|\leq |\mathbf{V}_b|+|\mathbf{V}_c|,\\
&\big||\mathbf{V}_c|-|\mathbf{V}_a|\big| \leq|\mathbf{V}_{ca}|\leq |\mathbf{V}_c|+|\mathbf{V}_a|.
\end{align*}

Denote $|\mathbf{V}_{abc}|=(|\mathbf{V}_a|,|\mathbf{V}_b|,|\mathbf{V}_c|)^{\top}$. Since $\epsilon\in(0,1)$, we also have
\begin{align}
 &\min_{V_{abc}\in\mathcal{V}_{in}} (\epsilon+2)|\mathbf{V}_{ab}| + (\epsilon-1)|\mathbf{V}_{bc}|+ (\epsilon-1)|\mathbf{V}_{ca}|\nonumber\\
 &\geq \min_{V_{abc}\in\mathcal{V}_{in}} (\epsilon+2)\big||\mathbf{V}_a|-|\mathbf{V}_b|\big| + (\epsilon-1)(|\mathbf{V}_b|+|\mathbf{V}_c|)\nonumber\\ &\hspace{2cm}+(\epsilon-1)(|\mathbf{V}_c|+|\mathbf{V}_a|)\nonumber\\
 &= \min_{Q|\mathbf{V}_{abc}|\leq q} (\epsilon+2)\big||\mathbf{V}_a|-|\mathbf{V}_b|\big| + (\epsilon-1)(|\mathbf{V}_b|+|\mathbf{V}_c|)\nonumber\\ &\hspace{2cm}+(\epsilon-1)(|\mathbf{V}_c|+|\mathbf{V}_a|)\label{eq: lvurconserv}
\end{align}
where $Q=[\mathbb{I}_3,-\mathbb{I}_3]^{\top}$ and 
\begin{align*}  
q=&(r_a+\|C_a\|,r_b+\|C_b\|,r_c+\|C_c\|, \min\{r_a-\|C_a\|,0\},\\
&\min\{r_b-\|C_b\|,0\},\min\{r_c-\|C_c\|,0\})^{\top}.
\end{align*}

The last equality in \eqref{eq: lvurconserv} is true since $V_{abc}$ is independent in each phase in $\mathcal{V}_{in}$, Similar to \eqref{eq: pvurrob2}, only the upper and lower bounds of $|\mathbf{V}_{abc}|$ are taking effect. It can be seen that \eqref{eq: lvurconserv} is convex and the feasible set $Q|\mathbf{V}_{abc}|\leq q$ only has eight extreme points (combinations of upper and lower bounds for $|\mathbf{V}_{abc}|$. Hence, evaluating the extreme points and finding the minimum efficiently solves \eqref{eq: lvurconserv}. Linear program duality \cite{Boyd_Convex_2004} can also be used to handle \eqref{eq: lvurconserv}. Below, we directly give the duality-based safe approximation to \eqref{eq: lvurrob} by introducing the dual variable $\lambda$:
\begin{align*}
-\hat{q}^{\top}\lambda\geq 0,\ \hat{Q}^{\top}\lambda+c=0,\ \lambda\geq 0,
\end{align*}
where $\hat{q}=(q^{\top},0,0)^{\top}$, $c=(\epsilon-1,\epsilon-1,2\epsilon-2,\epsilon+2)^{\top}$, and 
$$\hat{Q}=\begin{bmatrix}\begin{array}{cll|l}
\multicolumn{3}{c|}{\multirow{3}{*}{$Q$}} & \multirow{3}{*}{$\hphantom{-}\mathbf{0}$} \\
\multicolumn{3}{c|}{}                   &                    \\
\multicolumn{3}{c|}{}                   &                    \\ \hline
\multicolumn{1}{l}{\hphantom{-}1}     & -1    & 0    & -1                  \\
\multicolumn{1}{l}{-1}     & \hphantom{-}1    & 0    & -1                        
\end{array}\end{bmatrix}.$$

\subsection{Voltage Unbalance Factor (VUF) Definition}
References \cite{iec} and \cite{IEEEstd2} give the following voltage unbalance factor ($VUF$) definitions based on the magnitudes of \mbox{negative-,} \mbox{positive-,} and zero-sequence voltages, $\mathbf{V}_n$, $\mathbf{V}_p$, and $\mathbf{V}_0$, respectively:
\begin{subequations}
\begin{align}
    VUF_n=|\mathbf{V}_n|/|\mathbf{V}_p|,\\
    VUF_0=|\mathbf{V}_0|/|\mathbf{V}_p|,
\end{align}
\end{subequations}
where
\begin{subequations}\label{eq: seq}
\begin{align}
    &\mathbf{V}_p=(\mathbf{V}_a+\alpha\mathbf{V}_b+\alpha^2\mathbf{V}_c)/{3},\\
    &\mathbf{V}_n=(\mathbf{V}_a+\alpha^2\mathbf{V}_b+\alpha\mathbf{V}_c)/{3},\\
    &\mathbf{V}_0=(\mathbf{V}_a+\mathbf{V}_b+\mathbf{V}_c)/{3},
\end{align}
\end{subequations}
and $\alpha=1\angle 120$. With the tolerance $\epsilon\in(0,1)$, we equivalently transform the voltage balance requirements into quadratic inequality constraints:
\begin{subequations}\label{eq: vuf_complex}
\begin{align}
    &|\mathbf{V}_n|/|\mathbf{V}_p|\leq\epsilon
\quad \Leftrightarrow  \quad \mathbf{V}_n\bar{\mathbf{V}}_n-\epsilon^2\mathbf{V}_p\bar{\mathbf{V}}_p\leq 0,\\
 &|\mathbf{V}_0|/|\mathbf{V}_p|\leq\epsilon
\hspace*{0.85pt}\quad \Leftrightarrow  \quad \mathbf{V}_0\bar{\mathbf{V}}_0-\epsilon^2\mathbf{V}_p\bar{\mathbf{V}}_p\leq 0.
\end{align}
   \end{subequations}
   
Next, to ensure the power flow solutions are balanced, we obtain the following constraints
\begin{subequations}\label{eq: vufrob_complex}
\begin{align}
 & \left\{\max_{\mathbf{V}_{abc}\in\mathbf{U}_{in}} \mathbf{V}_n\bar{\mathbf{V}}_n-\epsilon^2\mathbf{V}_p\bar{\mathbf{V}}_p\right\}\leq 0,\label{eq: vufrob_c1} \\
 & \left\{\max_{\mathbf{V}_{abc}\in\mathbf{U}_{in}} \mathbf{V}_0\bar{\mathbf{V}}_0-\epsilon^2\mathbf{V}_p\bar{\mathbf{V}}_p\right\}\leq 0.\label{eq: vufrob_c2} 
\end{align}
  \end{subequations}
  
A direct way to safely approximate \eqref{eq: vufrob_complex} is using approximation by bound. For example, \eqref{eq: vufrob_c1} is implied by 
\begin{equation}
\left\{\max_{\mathbf{V}_{abc}\in\mathbf{U}_{in}}9 \mathbf{V}_n\bar{\mathbf{V}}_n-\epsilon^2\min_{\mathbf{V}_{abc}\in\mathbf{U}_{in}}9\mathbf{V}_p\bar{\mathbf{V}}_p\right\}\leq 0\label{eq: vufrob_complex_conserv}
\end{equation}
where scaling helps eliminate $1/3$ in \eqref{eq: seq}. Further, we have
\begin{align*}
  \max_{\mathbf{V}_{abc}\in\mathbf{U}_{in}} 9\mathbf{V}_n\bar{\mathbf{V}}_n \leq (|\mathbf{C}_a+\alpha^2\mathbf{C}_b+\alpha\mathbf{C}_c|+r_a+r_b+r_c)^2  
\end{align*}
and the inequality is tight when $\mathbf{V}_a-\mathbf{C}_a$, $\alpha^2(\mathbf{V}_b-\mathbf{C}_b)$, and $\alpha(\mathbf{V}_c-\mathbf{C}_c)$ share the same angle as $\mathbf{C}_a+\alpha^2\mathbf{C}_b+\alpha\mathbf{C}_c$.
Similarly, we get 
\begin{align*}
  & \min_{\mathbf{V}_{abc}\in\mathbf{U}_{in}} 9\mathbf{V}_p\bar{\mathbf{V}}_p \\ & \qquad =  (\max\{|\mathbf{C}_a+\alpha\mathbf{C}_b+\alpha^2\mathbf{C}_c|-r_a-r_b-r_c,0\})^2.  
\end{align*}

Hence, \eqref{eq: vufrob_complex_conserv} is equivalent to 
\begin{align}
    &(|\mathbf{C}_a+\alpha^2\mathbf{C}_b+\alpha\mathbf{C}_c|+r_a+r_b+r_c)^2\nonumber\\
    &\leq\ \epsilon^2(\max\{|\mathbf{C}_a+\alpha\mathbf{C}_b+\alpha^2\mathbf{C}_c|-r_a-r_b-r_c,0\})^2
\end{align}
and \eqref{eq: vufrob_c2} can be handled similarly. In addition to the approximation by bound, we give other approximation techniques by further transforming \eqref{eq: vufrob_complex} into the real domain using $V_{abc}$:
\begin{subequations}
\begin{align}
 & \left\{\max_{V_{abc}\in\mathcal{V}_{in}}V_{abc}^{\top}(A_n-\epsilon^2A_p)V_{abc}\right\}\leq 0,\label{eq: vufrob} \\
 & \left\{\max_{V_{abc}\in\mathcal{V}_{in}}V_{abc}^{\top}(A_0-\epsilon^2A_p)V_{abc}\right\}\leq 0.\label{eq: vufrob2} 
\end{align}
  \end{subequations}
where $V_{abc}^{\top}A_nV_{abc}=9\mathbf{V}_n\bar{\mathbf{V}}_n$, $V_{abc}^{\top}A_0V_{abc}=9\mathbf{V}_0\bar{\mathbf{V}}_0$, and $V_{abc}^{\top}A_pV_{abc}=9\mathbf{V}_p\bar{\mathbf{V}}_p$. Matrices $A_n\in\mathbb{R}^{6\times6}$, $A_p\in\mathbb{R}^{6\times6}$, and $A_0\in\mathbb{R}^{6\times6}$ can be calculated from \eqref{eq: seq} and have the following structure with off-diagonal matrices $B_n\in\mathbb{R}^{2\times2}$, $B_0\in\mathbb{R}^{2\times2}$, and $B_p\in\mathbb{R}^{2\times2}$
\begin{subequations}
\begin{align}
&A_n=\begin{bmatrix} \mathbb{I}_2 &B_n& B_n^{\top}\\B_n^{\top} & \mathbb{I}_2 & B_n\\
B_n&B_n^{\top}&\mathbb{I}_2
    \end{bmatrix},\ B_n=\begin{bmatrix}\cos(240) & -\sin(240)\\\sin(240) & \cos(240)\end{bmatrix},\\
&A_0=\begin{bmatrix} \mathbb{I}_2 &B_0& B_0^{\top}\\B_0^{\top} & \mathbb{I}_2 & B_0\\
B_0&B_0^{\top}&\mathbb{I}_2
    \end{bmatrix},\ B_0=\begin{bmatrix}\cos(0) & -\sin(0)\\\sin(0) & \cos(0)\end{bmatrix},\\
&A_p=\begin{bmatrix} \mathbb{I}_2 &B_p& B_p^{\top}\\B_p^{\top} & \mathbb{I}_2 & B_p\\
B_p&B_p^{\top}&\mathbb{I}_2
    \end{bmatrix},\ B_p=\begin{bmatrix}\cos(120) & -\sin(120)\\\sin(120) & \cos(120)\end{bmatrix}.
\end{align}
\end{subequations}


Both $A_n$ and $A_p$ are rank-two matrices and all four corresponding eigenvectors are orthogonal to each other. Hence, the matrix $A_n-\epsilon^2A_p$ is indefinite with rank four and the left-hand side (LHS) of \eqref{eq: vufrob} is a nonconvex quadratically constrained quadratic program (QCQP) with multiple constraints. A similar conclusion holds for the LHS of \eqref{eq: vufrob2}. General non-convex QCQPs are NP-hard to solve. 

To effectively approximate the QCQP or its solution, we first give the following lemma that provides a necessary condition on the location of the optimal solutions. For the rest of the paper, we use \eqref{eq: vufrob} and $VUF_n$ as an example since \eqref{eq: vufrob2} and $VUF_0$ can be similarly handled with exactly the same theoretical properties. 
\begin{lemma}\label{lem: boundary}
If $V^*_{abc}=(V_a^*,V_b^*,V_c^*)^{\top}$ is optimal for \eqref{eq: vufrob}, then $$V_a^*\in\partial\mathcal{V}^{a}_{in},\ V_b^*\in\partial\mathcal{V}^{b}_{in},\ V_c^*\in\partial\mathcal{V}^{c}_{in}.$$

\end{lemma}
\begin{IEEEproof}
We prove by contradiction. First, we assume that $V_a^*\notin\partial\mathcal{V}^{a}_{in}$ and define a correponding vector $\Delta_V=\hat{\alpha}(1,0,0,0,0,0)^{\top}\in\mathbb{R}^6$ with scalar $\hat{\alpha}$. Then, we conclude that there exists $\delta>0$ such that $(V_{abc}^*+\Delta_V)\in\mathcal{V}_{in}$ for all $\{\hat{\alpha}\in\mathbb{R}:\ |\hat{\alpha}|<\delta\}$ since $V_a^*\notin\partial\mathcal{V}^{a}_{in}$. Next, we compare the optimal objective with the objective under $(V_{abc}^*+\Delta_V)$:
\begin{align}
   &(V^*_{abc}+\Delta_V)^{\top}(A_n-\epsilon^2A_p)(V^*_{abc}+\Delta_V)\nonumber\\
   &\hspace{4cm} -(V^*_{abc})^{\top}(A_n-\epsilon^2A_p)V^*_{abc}\nonumber\\
   =&\ \hat{\alpha}^2(\Delta_V^{\top}(A_n-\epsilon^2A_p)\Delta_V)+\hat{\alpha}(2\Delta_V^{\top}(A_n-\epsilon^2A_p)V^*_{abc})\nonumber\\
   =&\ \hat{\alpha}^2(1-\epsilon^2)+\hat{\alpha}(2\Delta_V^{\top}(A_n-\epsilon^2A_p)V^*_{abc})=f(\hat{\alpha}).
\end{align}

Since $\epsilon<1$, we have $1-\epsilon^2>0$ and $f(\hat{\alpha})$ is a convex quadratic function of $\hat{\alpha}$. When $\hat{\alpha}=0$, we have $f(0)=0$ and hence we must also have $\max\{f(\frac{\delta}{2}),f(-\frac{\delta}{2})\}>0$. In other words, we can improve the optimal value of \eqref{eq: vufrob} by choosing either $\hat{\alpha}=\frac{\delta}{2}$ or $-\frac{\delta}{2}$ and constructing a new solution $(V_{abc}^*+\Delta_V)\in\mathcal{V}_{in}$. Hence, this is contradictory with $V^*_{abc}$ being optimal. Similar discussions apply to cases when $V_b^*\notin\partial\mathcal{V}^{b}_{in}$ and $V_c^*\notin\partial\mathcal{V}^{c}_{in}$ and the proof is complete.  
\end{IEEEproof}

Using Lemma \ref{lem: boundary}, we know that \eqref{eq: vufrob} is equivalent to
\begin{equation}
 \left\{\max_{V_{abc}\in\partial\mathcal{V}^a_{in}\times\partial\mathcal{V}^b_{in}\times\partial\mathcal{V}^c_{in}}V_{abc}^{\top}(A_n-\epsilon^2A_p)V_{abc}\right\}\leq 0.\label{eq: vufrob_bound}  
\end{equation}
We develop two approaches to approximate \eqref{eq: vufrob_bound} or its solution. 


\subsubsection{Polytope Approximation}
First, we model three polytopes $\mathcal{P}^a\in\mathbb{R}^2$, $\mathcal{P}^b\in\mathbb{R}^2$, and $\mathcal{P}^c\in\mathbb{R}^2$ such that
\begin{align}
    \mathcal{V}_{in}^a\subset\mathcal{P}^a,\ \mathcal{V}_{in}^b\subset\mathcal{P}^b,\ \mathcal{V}_{in}^c\subset\mathcal{P}^c.
\end{align} 

Denote the finite set of the extreme points of $\mathcal{P}^a$, $\mathcal{P}^b$, and $\mathcal{P}^c$ as $\mathcal{E}^a$, $\mathcal{E}^b$, and $\mathcal{E}^c$, respectively. We next present a theorem that provides a necessary condition on the location of the optimal solution for the maximization problem in
\begin{align}
 \left\{\max_{V_{abc}\in\mathcal{P}^a\times\mathcal{P}^b\times\mathcal{P}^c}V_{abc}^{\top}(A_n-\epsilon^2A_p)V_{abc}\right\}\leq 0.\label{eq: vufrobpoly}  
\end{align}
It is easy to see that \eqref{eq: vufrobpoly} is a safe approximation of \eqref{eq: vufrob} and \eqref{eq: vufrob_bound} with a larger optimal value.

\begin{thm} \label{thm: bound}
If $V^*_{abc}=(V^*_a,V^*_b,V^*_c)^{\top}$ is optimal for the maximization problem in \eqref{eq: vufrobpoly}, then
\begin{align}
    V^*_a\in\mathcal{E}^a,\ V^*_b\in\mathcal{E}^b,\ V^*_c\in\mathcal{E}^c.
\end{align}
\end{thm}  
\begin{IEEEproof}
First, we claim that 
\begin{align}
    V^*_a\in\partial\mathcal{P}^a,\ V^*_b\in\partial\mathcal{P}^b,\ V^*_c\in\partial\mathcal{P}^c
\end{align}
whose proof is similar to the one of Lemma~\ref{lem: boundary}.

Next, we show the theorem by contradiction. Since $V^*_a\in\partial\mathcal{P}^a$, then $V^*_a\in H^a$ where $H^a$ is one of the hyperplanes defining $\partial\mathcal{P}^a$. Define $H^a$ as $\{x\in\mathbb{R}^2:\ h^{\top}x=\tilde{h}\}$, then $h^{\top}V^*_a=\tilde{h}$. If we assume $V_a^*\notin\mathcal{E}^a$, then there exists $\delta>0$ and a direction $\{g\in\mathbb{R}^2: \|g\|=1,\ g^{\top}h=0\}$ such that \mbox{$(V^*_a+\hat{\alpha} g)\in\partial\mathcal{P}^a$} for all $\{\hat{\alpha}\in\mathbb{R}:\ |\hat{\alpha}|<\delta\}$. Next, we compare the optimal objective with the objective under $(V^*_{abc}+\Delta_V)$ where $\Delta_V=\hat{\alpha}(g^{\top},0,0,0,0)^{\top}\in\mathbb{R}^6$ and get
\begin{align}
   &(V^*_{abc}+\Delta_V)^{\top}(A_n-\epsilon^2A_p)(V^*_{abc}+\Delta_V) \nonumber\\
   &\hspace{4cm} -(V^*_{abc})^{\top}(A_n-\epsilon^2A_p)V^*_{abc}\nonumber\\
   =&\ \hat{\alpha}^2(\Delta_V^{\top}(A_n-\epsilon^2A_p)\Delta_V)+\hat{\alpha}(2\Delta_V^{\top}(A_n-\epsilon^2A_p)V^*_{abc}) \nonumber\\
   =&\ \hat{\alpha}^2(1-\epsilon^2)g^{\top}\mathbb{I}_2g+\hat{\alpha}(2\Delta_V^{\top}(A_n-\epsilon^2A_p)V^*_{abc})\nonumber\\
   =&\ \hat{\alpha}^2(1-\epsilon^2)+\hat{\alpha}(2\Delta_V^{\top}(A_n-\epsilon^2A_p)V^*_{abc})=f(\hat{\alpha}),
\end{align}
which is a convex quadratic function on $\hat{\alpha}$ with $f(0)=0$ since $\epsilon<1$. Then, similar to Lemma \ref{lem: boundary}, we conclude that we can improve the optimal value of \eqref{eq: vufrobpoly} by using a new feasible solution $(V^*_{abc}+\Delta_V)$ with $\hat{\alpha}=\frac{\delta}{2}$ or $-\frac{\delta}{2}$. This contradicts the optimality of $V^*_{abc}$.
Similar discussions are applicable to cases when $V_b^*\notin\mathcal{E}^{b}$ and $V_c^*\notin\mathcal{E}^{c}$ and the proof is complete.
\end{IEEEproof}

Now, we equivalently reformulate \eqref{eq: vufrobpoly} as
\begin{equation}
 \left\{\max_{V_{abc}\in\mathcal{E}^a\times\mathcal{E}^b\times\mathcal{E}^c}V_{abc}^{\top}(A_n-\epsilon^2A_p)V_{abc}\right\}\leq 0\label{eq: vufrobextre}  
\end{equation}
and solving an optimization problem \eqref{eq: vufrobpoly} becomes an evaluation problem on the set of extreme points. 


There are many ways to find $\mathcal{P}^a$, $\mathcal{P}^b$, and $\mathcal{P}^c$. Here, we use a special polytope to analyze the optimality gap of the approximation. Since each polytope is in dimension $2$, we propose to use the circumscribed regular polygon of the disk (CRP). For a unit closed disk $\bar{D}(1)$, the extreme points of a CRP with $2m$ ($m\geq2$) sides are as follows 
$$\left\{\frac{1}{\cos(\frac{\pi}{2m})}\begin{bmatrix}\cos(\phi)\\ \sin(\phi) \end{bmatrix}:\ \phi=\frac{(2k-1)\pi}{2m},\ k=1,2,...,2m\right\}.$$

Note that while a CRP can have a phase shift, we do not consider this here for the sake of simplicity. In combination with $C_{abc}$ and $r_{abc}$, we can easily find $\mathcal{E}^{a,2m}$, $\mathcal{E}^{b,2m}$, and $\mathcal{E}^{c,2m}$. We add $2m$ in the notations to denote the dimension of the CRP. By defining a general function $E^{2m}:\mathbb{R}^2\times\mathbb{R}\to \mathbb{R}^{2m}$, then $\mathcal{E}^a=E^{2m}(C_a,r_a)$ can be represented as
\begin{align*}
\Bigg\{C_a+\frac{r_a}{\cos(\frac{\pi}{2m})}&\begin{bmatrix}\cos(\phi)\\ \sin(\phi) \end{bmatrix}:\ \\
&\phi=\frac{(2k-1)\pi}{2m},\ k_a=1,2,...,2m\Bigg\}.    
\end{align*}

We next show how the optimality gap between \eqref{eq: vufrob_bound} and \eqref{eq: vufrobextre} is affected by the dimension $m$.  

\begin{coro} Denoting the optimal values of \eqref{eq: vufrob_bound} and \eqref{eq: vufrobextre} as $F^*_{b}$ and $F^*_{e}$, respectively, we have
\begin{equation}
  |F^*_{e}-F^*_{b}| \leq |F^*_{e}-F^*_{i}|,  \label{eq: gap}
\end{equation}
where $F^*_{i}$ is the optimal solution of the following problem
\begin{equation}
    \max_{V_{abc}\in\hat{\mathcal{E}}^a\times\hat{\mathcal{E}}^b\times\hat{\mathcal{E}}^c}V_{abc}^{\top}(A_n-\epsilon^2A_p)V_{abc}, \label{eq: vufrobextre2}
\end{equation}
in which 
\begin{align*}
    &\hat{\mathcal{E}}^a=E^{2m}(C_a,r_a\cos(\frac{\pi}{2m})),\\
    &\hat{\mathcal{E}}^b=E^{2m}(C_b,r_b\cos(\frac{\pi}{2m})),\\
    &\hat{\mathcal{E}}^c=E^{2m}(C_c,r_c\cos(\frac{\pi}{2m})).
\end{align*}

We also have 
\begin{equation}
    \lim_{m\to+\infty}|F^*_{e}-F^*_{i}|=0. \label{eq: gapzero}
\end{equation}
\end{coro}
\begin{IEEEproof} We prove \eqref{eq: gap} by demonstrating the relationship
$$F^*_{e}\geq F^*_{b} \geq F^*_{i}.$$

The first inequality results from the fact that \eqref{eq: vufrobextre} is a safe approximation of \eqref{eq: vufrob_bound}. The second inequality is true because $\hat{\mathcal{E}}^a\subset\partial\mathcal{V}_{in}^a$ (same for phases $b$ and $c$). Hence, \eqref{eq: vufrob_bound} is more conservative than \eqref{eq: vufrobextre2}. 

Next, we prove \eqref{eq: gapzero}. Given any $m\geq 2$, we have the following inequality
\begin{equation}
    |F^*_{e}-F^*_{i}|\leq \max_{
         k_a\in\mathcal{K},
         k_b\in\mathcal{K},
         k_c\in\mathcal{K}} |J(V^e_{abc})-J(V^i_{abc})|,\label{eq: ineq1}
\end{equation}
where $J(V_{abc})=V_{abc}^{\top}(A_n-\epsilon^2A_p)V_{abc}$. $V^e_{abc}$ and $V^i_{abc}$ are a corresponding pair in $\mathcal{E}^a\times\mathcal{E}^b\times\mathcal{E}^c$ and $\hat{\mathcal{E}}^a\times \hat{\mathcal{E}}^b\times \hat{\mathcal{E}}^c$ with the same $k_a$, $k_b$, and $k_c$. $\mathcal{K}$ denotes the integer set on $[1,2m]$. We can see \eqref{eq: ineq1} is valid if we substitute the $(k_a,k_b,k_c)^{\top}$, which is optimal for \eqref{eq: vufrobextre}, into the right-hand side of \eqref{eq: ineq1}. Further, we have
\begin{equation}
  \|V^e_{abc}-V^i_{abc}\|=\left(\frac{1}{\cos(\frac{\pi}{2m})}-1\right)\|r_{abc}\|,\ \forall (k_a,k_b,k_c)^{\top}\in\mathcal{K}^3 \label{eq: ineq3}
\end{equation}
since $V^e_{abc}$ and $V^i_{abc}$ are a corresponding pair. Meanwhile, $J(V_{abc})$ is continuously differentiable and hence Lipschitz on compact set $\bar{D}(C_a,\sqrt{2}r_a)\times \bar{D}(C_b,\sqrt{2}r_b)\times \bar{D}(C_c,\sqrt{2}r_c)$. We choose this compact set since it contains all the feasible sets of \eqref{eq: vufrobextre} and \eqref{eq: vufrobextre2} for all $m\geq 2$. Denote the Lipschitz constant as $L$.  For all $m\geq 2$, we have
\begin{align}
&\max_{
         k_a\in\mathcal{K},
         k_b\in\mathcal{K},
         k_c\in\mathcal{K}} |J(V^e_{abc})-J(V^i_{abc})|\nonumber\\
         & \qquad \leq \max_{
         k_a\in\mathcal{K},
         k_b\in\mathcal{K},
         k_c\in\mathcal{K}} L\|V^e_{abc}-V^i_{abc}\|.\label{eq: ineq2}
\end{align}

Combining \eqref{eq: ineq1}, \eqref{eq: ineq3}, and \eqref{eq: ineq2}, we have
\begin{align*}
 &\lim_{m\to+\infty}|F^*_{e}-F^*_{i}|\leq \lim_{m\to+\infty}\max_{
         k_a\in\mathcal{K},
         k_b\in\mathcal{K},
         k_c\in\mathcal{K}} L\|V^e_{abc}-V^i_{abc}\|\\
         &= \lim_{m\to+\infty} (\frac{1}{\cos(\frac{\pi}{2m})}-1)L\|r_{abc}\|=0.    
\end{align*}

The last equality holds since $\lim_{m\to+\infty}\cos(\frac{\pi}{2m})=1$. Further, since $|F^*_{e}-F^*_{i}|$ is non-negative, based on squeeze theorem, the proof is complete. 
\end{IEEEproof}

This result tells us that as we increase $m$ (i.e., the number of sides of the CRP), the safe approximation \eqref{eq: vufrobextre} asymptotically converges to the true optimal value of \eqref{eq: vufrob} or \eqref{eq: vufrob_bound}. 


\subsubsection{Semidefinite and Lagrangian Relaxation}
Other conventional techniques \cite{qcqp, Boyd_Convex_2004} for general QCQP problems use semidefinite relaxation (SDR) or Lagrangian relaxation (LGR). The SDR of \eqref{eq: vufrob_bound} (shown below) is derived by lifting the vector space of the variable $V_{abc}$ to the matrix space $W_{abc}\in\mathbb{R}^{6\times 6}$ and relaxing the rank-one constraints from $W_{abc}=V_{abc}V_{abc}^{\top}$ to get a convex constraint $W_{abc}\succeq V_{abc}V_{abc}^{\top}$ and the following semidefinite programming (SDP) problem. Since the original problem in \eqref{eq: vufrob_bound} is maximization, both LGR and SDR give higher optimal values and hence a safe approximation to \eqref{eq: vufrob_bound}:
\begin{align*}
    \mbox{(SDR)}\ \max\ &\Tr((A_n-\epsilon^2A_p)W_{abc})\\
    \mbox{s.t.}\ &W_{abc,11}+W_{abc,22}-2C_a^{\top}V_a+\|C_a\|^2=r_a^2,\\
                 &W_{abc,33}+W_{abc,44}-2C_b^{\top}V_b+\|C_b\|^2=r_b^2,\\
                 &W_{abc,55}+W_{abc,66}-2C_c^{\top}V_a+\|C_c\|^2=r_c^2,\\
                 &W_{abc}\succeq V_{abc}V_{abc}^{\top}
\end{align*}

LGR uses Lagrangian duality to derive an SDP-based reformulation as follows 
\begin{align*}
    \mbox{(LGR)}\ \min_{\mu\in\mathbb{R}^3}\ &\gamma\\
    \mbox{s.t.}\ &Y=\begin{bmatrix} Q(\mu) & q(\mu)\\ q(\mu)^{\top} & r(\mu)\end{bmatrix},\\
    &Y\succeq 0
\end{align*}
where 
\begin{align*}
    &Q(\mu)=-(A_n-\epsilon^2A_p+\text{blkdiag}(\mu_1\mathbb{I}_2,\mu_2\mathbb{I}_2,\mu_3\mathbb{I}_2)),\\
    &q(\mu)=(-\mu_1C_a^{\top}, -\mu_2C_b^{\top}, -\mu_3C_c^{\top})^{\top},\\
    &r(\mu)=\gamma-\mu_1(\|C_a\|^2-r_a^2)\\
    &\hspace{3cm}-\mu_2(\|C_b\|^2-r_b^2)-\mu_3(\|C_c\|^2-r_c^2).
\end{align*}

Then, \eqref{eq: vufrob_bound} can be safely approximated as
\begin{align}
    \gamma \leq 0,\ \mu\in\mathbb{R}^3,\ \begin{bmatrix} Q(\mu) & q(\mu)\\ q(\mu)^{\top} & r(\mu)\end{bmatrix}\succeq 0. \label{eq: lgrfinal}
\end{align}

Existing work \cite{qcqp} shows that SDR and LGR are dual to each other. Strong duality also holds here as both SDR and LGR are strictly feasible (i.e., there exists positive definite matrix solutions).\footnote{In SDR, we can select $V_{abc}=0$ and pick $W_{abc}$ to be diagonal with strictly positive elements. In LGR, we can choose any $\mu < -\lambda_{\max}(A_n-\epsilon^2A_p) $, then $Q(\mu)$ is positive definite. Then, based on the Schur complement, we can always choose $\gamma$ large enough that $r(\mu)-q(\mu)^{\top}P(\mu)^{-1}q(\mu)>0$.} Since \eqref{eq: vufrob_bound} is nonconvex, there is a gap between SDR or LGR with the true optimal solutions. Next, we give conditions on $C_{abc}$ and $r_{abc}$ such that SDR and LGR have the same optimal value to \eqref{eq: vufrob_bound}.  

We start from LGR and show a sufficient condition such that strong duality holds between LGR and \eqref{eq: vufrob_bound}. The conditions can also be efficiently evaluated by solving three small convex QCQPs. For concise derivation, we define three sets:
\begin{align*}
&\bar{D}_{abc}=\bar{D}(r_a)\times\bar{D}(r_b)\times\bar{D}(r_c),\\
&D_{abc}=D(r_a)\times D(r_b) \times D(r_c),\\
&\partial\bar{D}_{abc}=\partial\bar{D}(r_a)\times\partial\bar{D}(r_b)\times\partial\bar{D}(r_c).
\end{align*}

\begin{thm}\label{thm: suff1}
If conditions 
\begin{subequations}
\begin{align}
    &f_{i,i+1}\notin[-r_i,r_i]\times[-r_{i+1},r_{i+1}],\ i=1,3,5,\\
    &\left\{\min_{Y_a\in\bar{D}_{abc}}\|(2BY_a+f)_{1,2}\|^2\right\}\geq 4(2+\epsilon^2)^2 r_a^2,\\
    &\left\{\min_{Y_b\in\bar{D}_{abc}}\|(2BY_b+f)_{3,4}\|^2\right\}\geq 4(2+\epsilon^2)^2 r_b^2,\\
    &\left\{\min_{Y_c\in\bar{D}_{abc}}\|(2BY_c+f)_{5,6}\|^2\right\}\geq 4(2+\epsilon^2)^2 r_c^2,
\end{align}
\end{subequations}
where 
\begin{align*}
&B=(\epsilon^2A_p-A_n)-\lambda_{\min}(\epsilon^2A_p-A_n)\mathbb{I}_6,\\
&f=2(\epsilon^2A_p-A_n)C_{abc},\\
&r_i=2(r_a\|B_i^{1,2}\|+r_b\|B_i^{3,4}\|+r_c\|B_i^{5,6}\|),\ i=1,..,6, \end{align*}%
are satisfied by certain $C_{abc}$ and $r_{abc}$, then strong duality holds between LGR and \eqref{eq: vufrob_bound}.
\end{thm}
\begin{IEEEproof}
We prove strong duality for the following problem:
\begin{align}
\min_{V_{abc}\in\partial\mathcal{V}^a_{in}\times\partial\mathcal{V}^b_{in}\times\partial\mathcal{V}^c_{in}}V_{abc}^{\top}(\epsilon^2A_p-A_n)V_{abc}.\label{eq: dual}
\end{align}

If strong duality holds for \eqref{eq: dual}, strong duality also holds for \eqref{eq: vufrob_bound} since their duals always have opposite optimal values. 
Define $Y_{abc}=(Y_a,Y_b,Y_c)^{\top}=V_{abc}-C_{abc}$ and $A=\epsilon^2A_p-A_n$. We rewrite \eqref{eq: dual1} as
\begin{align}
\min_{Y_{abc}\in\partial\bar{D}_{abc}}(Y_{abc}+C_{abc})^{\top}A\,(Y_{abc}+C_{abc}),\nonumber
\end{align}
where
\begin{align*}
    &(Y_{abc}+C_{abc})^{\top}A\,(Y_{abc}+C_{abc})\\
    =&Y_{abc}^{\top}(A-\lambda\mathbb{I}_6)Y_{abc} + 2C_{abc}^{\top}AY_{abc}+C_{abc}^{\top}AC_{abc} + \lambda Y_{abc}^{\top}Y_{abc}.
\end{align*}

Since $Y_{abc}\in\partial\bar{D}_{abc}$, we ignore the last two terms (constant-valued) and pick $\lambda=\lambda_{\min}(A)=-3$ such that $B=(A-\lambda\mathbb{I}_6)\succeq 0$. Define $f=2AC_{abc}$. We then obtain 
\begin{align}
\min_{Y_{abc}\in\partial\bar{D}_{abc}}Y_{abc}^{\top}B\,Y_{abc}+f^{\top}Y_{abc}.\label{eq: compact}
\end{align}

A sufficient condition ensuring that strong duality holds for \eqref{eq: compact} is when the following problem
\begin{equation}
\min_{Y_{abc}\in\bar{D}_{abc}}Y_{abc}^{\top}B\,Y_{abc}+f^{\top}Y_{abc}\label{eq: compact2}    
\end{equation}
has its optimal solution on $\partial\bar{D}_{abc}$ under certain requirements on $C_{abc}$ and $r_{abc}$. We leave the proof of this to the following lemma. Now, we assume $Y^*=(Y^*_a,Y^*_b,Y^*_c)^{\top}$ to be any point in $\bar{D}_{abc}$ and $Y_a^*\notin \partial\bar{D}(r_a)$ without loss of generality since similar discussions can apply to $Y_b^*$ and $Y_c^*$. The sufficient condition for strong duality is satisfied if there exists a direction $\{g_a\in \mathbb{R}^2:\ \|g_a\|=1\}$ in the space of $Y_a$ such that the objective can be improved when $Y_a^*$ moves to $\partial\bar{D}(r_a)$. Denote $\Delta_Y=d(g_a,0,0,0,0)^{\top}$ where $d\geq0$ denotes the distance moving along the direction $\Delta_Y$. We compare the objective at $Y^*+\Delta_Y$ and $Y^*$ and get 
\begin{align}
   &(Y^*+\Delta_Y)^{\top}B\,(Y^*+\Delta_Y)+f^{\top}(Y^*+\Delta_Y)\nonumber\\
   &\hspace{3cm} -(Y^*)^{\top}B\,Y^*-f^{\top}Y^*\nonumber\\
   =&\ d^2(2+\epsilon^2)+d\,h_a^{\top}g_a=f(d,g_a)
\end{align}
where $h_a=(2BY^*+f)_{1,2}$. We cannot have $h_a=0$ since there does not exist a direction $g_a$ that improves the objective. Denote $(2BY^*)_i = 2B_i^{\top}Y^*=2B_i^{1,2}Y^*_a + 2B_i^{3,4}Y^*_b + 2B_i^{5,6}Y^*_c$. Since $Y^*$ is in a compact set $ \bar{D}_{abc}$, we can find the tight bound $|(2BY^*)_i|\leq r_i=2(r_a\|B_i^{1,2}\|+r_b\|B_i^{3,4}\|+r_c\|B_i^{5,6}\|)$.   
Further, if $f_{1,2}\notin[-r_1,r_1]\times[-r_2,r_2]$, then $h_a\neq 0$ for all possible $Y^*\in \bar{D}_{abc}$.


Next, since $f(d,g_a)$ is a convex quadratic function of $d$ with fixed second-order coefficient, we know the steepest descent direction is $-h_a/\|h_a\|$ (i.e., for any fixed $d$, $g_a=-h_a/\|h_a\|$ minimizes $f(d,g_a)$). Hence, we fix $g_a=-h_a/\|h_a\|$ without loss of optimality. The objective is improved if $f(d,-h_a/\|h_a\|)\leq 0$, which is equivalent to $d\in[0,\frac{\|h_a\|}{2+\epsilon^2}]$. Hence, if we require $\frac{\|h_a\|}{2+\epsilon^2}\geq 2r_a$, then there must exists a point on $\partial\bar{D}(r_a)$ with certain $d^*\in(0,\frac{\|h_a\|}{2+\epsilon^2})$ that improves the objective. Since $Y^*$ can be any point in $\bar{D}_{abc}$, our requirement becomes
\begin{align}
    &\left\{\min_{Y^*\in\bar{D}_{abc}}\|h_a\|\right\}\geq 2(2+\epsilon^2) r_a\nonumber\\
    \Leftrightarrow&\left\{\min_{Y^*\in\bar{D}_{abc}}h_a^{\top}h_a\right\}\geq 4(2+\epsilon^2)^2 r_a^2 \label{eq: strong1}.
\end{align}

The left-hand side of \eqref{eq: strong1} is a convex QCQP\footnote{Strong duality holds since Slater's condition is satisfied (i.e., pick $Y^*\in D_{abc}$). Hence, \eqref{eq: strong1} can also be equivalently transformed into SDP constraints.} and hence can be solved easily. Similar analyses for phases $b$ and $c$ complete the proof.
\end{IEEEproof}

\begin{lemma}\label{lem: supp}
For given $C_{abc}$ and $r_{abc}$, if \eqref{eq: compact2} has its optimal solution on $\partial\bar{D}_{abc}$, then \eqref{eq: compact} has strong duality.
\end{lemma}

\begin{IEEEproof}
For simplicity, generalize \eqref{eq: compact} as
$$\min f(x)\ \mbox{s.t. } g_i(x) = 0,\ i=1,...,m,$$
which is equivalent to 
$$v_1^*=\{\min f(x)\ \mbox{s.t. } g_i(x) \geq 0,\ g_i(x) \leq 0,\ i=1,...,m\}$$
with optimal value $v_1^*$. The associated Lagrangian dual is 
\begin{equation}
  h_1^*=\max_{\lambda\geq0,\mu\geq0}\inf_{x}f(x)+\Sigma^m_i\lambda_ig_i(x)-\Sigma^m_i\mu_ig_i(x)  \label{eq: dual1}
\end{equation}
with optimal value $h_1^*$. Similarly, \eqref{eq: compact2} can be represented as 
$$v_2^*=\{\min f(x)\ \mbox{s.t. } g_i(x) \leq 0\ i=1,...,m\}$$
with the following Lagrangian dual problem
\begin{equation}
   h_2^*=\max_{\lambda\geq0}\inf_{x}f(x)+\Sigma^m_i\lambda_ig_i(x) \label{eq: dual2}
\end{equation} 
and has optimal values $v_2^*$ ($h_2^*$ for \eqref{eq: dual2}).

Denote the optimal solution to \eqref{eq: dual2} as $x^*$ and $\lambda^*$. Observe that \eqref{eq: dual2} is a special case of \eqref{eq: dual1} with $\mu=0$. Hence, we have $h_1^*\geq h_2^*$. Since \eqref{eq: compact2} is a convex QCQP with nonempty interior, we also have strong duality such that $v_2^*=h_2^*$. Meanwhile, we also have $v_1^*=v_2^*$ (since \eqref{eq: compact2} has its optimal solution on $\partial\bar{D}_{abc}$) and $h_1^*\leq v_1^*$ (weak duality). We conclude $h_1^*= v_1^*$ and hence strong duality holds for \eqref{eq: compact}.   
\end{IEEEproof}

Theorem \ref{thm: suff1} implies that when $C_{abc}$ has relatively larger magnitudes than $r_{abc}$, there is a higher chance of having strong duality between LGR and \eqref{eq: vufrob_bound}. Next, we start with SDR and show that under certain conditions, even if $r_{abc}$ has large magnitudes, we have exactness (i.e., the same optimal values) between SDR and \eqref{eq: vufrob_bound}. First, we give a general result for a QCQP problem whose SDR is exact. 

\begin{thm}\label{thm: suff2}
Consider the following QCQP problem on $x_{abc}=(x_a,x_b,x_c)^{\top}\in\mathbb{R}^6$:  
\begin{align}
    \max_{x_{abc}\in\partial\bar{D}_{abc}} x_{abc}^{\top}Ax_{abc} \label{eq: homo1}
\end{align}
where $A$ has the following structure 
$$
A=\begin{bmatrix} \lambda\mathbb{I}_2 &B& B^{\top}\\B^{\top} & \lambda\mathbb{I}_2 & B\\
B&B^{\top}&\lambda\mathbb{I}_2
    \end{bmatrix},\ B=\begin{bmatrix}b_1 & -b_2\\b_2 & b_1\end{bmatrix}.$$ 
    
Then, the following SDP relaxation is exact for \eqref{eq: homo1}:
    \begin{align}
    \max\ &\Tr(AX_{abc})\label{eq: sdr_homo}\\
    \mbox{s.t.}\ &X_{abc,11}+X_{abc,22}=r_a^2,\nonumber\\
                 &X_{abc,33}+X_{abc,44}=r_b^2,\nonumber\\
                 &X_{abc,55}+X_{abc,66}=r_c^2,\nonumber\\
                 &X_{abc}\succeq 0.\nonumber
\end{align}
\end{thm}

\begin{IEEEproof}
We prove the exactness by rewriting \eqref{eq: homo1} into a homogeneous complex QCQP with the transformation $\mathbf{x}_{abc}=(x_{a,1}+\bm{j}x_{a,2},x_{b,1}+\bm{j}x_{b,2},x_{c,1}+\bm{j}x_{c,2})^{\top}$:
\begin{align}
    \eqref{eq: homo1}\ \Leftrightarrow\ \max\ &\mathbf{x}_{abc}^{\text{H}}\mathbf{A}\mathbf{x}_{abc}\label{eq: homo2}\\
    \mbox{s.t. }&\mathbf{x}_{abc}^{\text{H}}[\text{blkdiag}(1,\mymathbb{0}_2)]\mathbf{x}_{abc}=r_a^2,\nonumber\\
    &\mathbf{x}_{abc}^{\text{H}}[\text{blkdiag}(0,1,0)]\mathbf{x}_{abc}=r_b^2,\nonumber\\
    &\mathbf{x}_{abc}^{\text{H}}[\text{blkdiag}(\mymathbb{0}_2,1)]\mathbf{x}_{abc}=r_c^2,\nonumber
\end{align}
where $$\mathbf{A}=\begin{bmatrix} \lambda &b_1+\bm{j}b_2& b_1-\bm{j}b_2\\b_1-\bm{j}b_2 & \lambda& b_1+\bm{j}b_2\\
b_1+\bm{j}b_2&b_1-\bm{j}b_2&\lambda
    \end{bmatrix}.$$
    
Then, based on Theorem 3.2 in \cite{sdp_complex1,sdp_complex2}, the SDP relaxation of \eqref{eq: homo2} has a solution $\mathbf{X}_{abc}\in\mathcal{H}^3$ with $\mathrm{rank}(\mathbf{X}_{abc})\leq \floor{\sqrt{l}}$ where $l$ equals the number of equality constraints in \eqref{eq: homo2} (i.e., $l=3$).  
Hence, $\mathrm{rank}(\mathbf{X}_{abc})\leq 1$ and is thus exact for \eqref{eq: homo2}. Since \eqref{eq: sdr_homo} is a real-variable representation of the SDP relaxation of \eqref{eq: homo2}, the proof is complete.\footnote{Note that exactness is not guaranteed for any other QCQP formulations (e.g., \eqref{eq: vufrob_bound}) that are either non-homogeneous or involve more than three equality constraints in complex-variable representation.}
\end{IEEEproof}


To see how Theorem~\ref{thm: suff2} is relevant to \eqref{eq: vufrob_bound}, we reformulate \eqref{eq: vufrob_bound} by substituting $Y_{abc}=V_{abc}-C_{abc}$ and ignoring the constant to obtain
\begin{align}
    \max_{Y_{abc}\in\partial\bar{D}_{abc}}Y_{abc}^{\top}(A_n-\epsilon^2A_p)Y_{abc}+2((A_n-\epsilon^2A_p)C_{abc})^{\top}Y_{abc},\nonumber
\end{align}
which satisfies Theorem~\ref{thm: suff2} if $(A_n-\epsilon^2A_p)C_{abc}=0$. Since strong duality holds between LGR and SDR, Theorem~\ref{thm: suff2} also guarantess strong duality between LGR and \eqref{eq: vufrob_bound}. Further, we see that Theorems \ref{thm: suff1} and \ref{thm: suff2} are not directly comparable (i.e., it is not the case that one is always stronger than the other). Theorem \ref{thm: suff1} can handle cases when $C_{abc}$ has large magnitude relative to $r_{abc}$. However, Theorem \ref{thm: suff2} works when $C_{abc}$ has a small magnitude compared to $A_n-\epsilon^2A_p$ or belongs to $N(A_n-\epsilon^2A_p)$ (approximately $N(A_n)$ when $\epsilon$ is small). Hence, each condition has its own advantages. In practice, Theorem \ref{thm: suff1} is more applicable as the condition is better aligned with the characteristics of practical distribution networks. 

\subsection{Balancibility Condition}
Using the closed-form reformulations or approximations of the voltage balance requirements from the three different unbalance definitions, we can now obtain the full balancibility condition by connecting the $\mathcal{V}_L$ from the solvability condition \eqref{eq: invset2} to these voltage balance requirements. We briefly discussed this connection when we first introduced $C_{abc}$ and $r_{abc}$ in Section~\ref{sec:bal}. Here, we give an example of the complete balancibility condition using the polytope approximation under the $VUF_n$ definition. Without loss of generality, assume $i^*\in\mathcal{N}_L$ is the critical node with a voltage unbalance tolerance of $\epsilon$ and there is no degeneracy (i.e., $\min_{p\in\{a,b,c\}}\xi_{i^*,p}(\mathbf{S}_L)>0$). Then, the full balancibility condition is:
\begin{align}
    \eqref{eq:solv}\ \text{and}\ \left\{\max_{V_{abc}\in\mathcal{E}^a\times\mathcal{E}^b\times\mathcal{E}^c}V_{abc}^{\top}(A_n-\epsilon^2A_p)V_{abc}\right\}\leq 0
\end{align}
where 
$$\mathcal{E}^p=E^{2m}\left((\re(\mathbf{C}_p),\im(\mathbf{C}_p))^{\top},r_p\right),$$ $\mathbf{C}_p=(1-\eta_{i^*,p}(\bm{\sigma}_L))\mathbf{E}_{p}^{i^*}\mathbf{v}_{L,p}^{0,i^*}$, and $r_p=r|\mathbf{E}_{p}^{i^*}\mathbf{v}_{L,p}^{0,i^*}|\xi_{i^*,p}(\mathbf{S}_L)$ for all $p\in\{a,b,c\}$. Note that the solvability condition can be simultaneously combined with multiple voltage balance requirements. Particularly for $VUF$ definitions, since $C_{abc}$ and $r_{abc}$ can be determined if $\mathbf{S}_L$ is provided, the sufficient conditions in Theorems \ref{thm: suff1} and \ref{thm: suff2} can be evaluated to give better guidance regarding whether we should choose a polytope approximation or LGR.

All of the constraints in the balancibility condition are functions of $\mathbf{S}_L$. Hence, the balancibility condition defines a secure region of $\mathbf{S}_L$ such that a unique and balanced power flow solution is guaranteed to exist. There are multiple applications of the balancibility condition, such as directly incorporating the condition in a centralized problem. This condition then provides a feasible set on $\mathbf{S}_L$ that can directly replace the power flow equations to obtain associated robust voltage balance guarantees for the power flow solutions under uncertainty. Another application uses the balancibility condition in an iterative/decentralized algorithm such that each step provides an instance of $\mathbf{S}_L$. In this case, the balancibility condition can be used as an effective evaluation tool on $\mathbf{S}_L$ to provides quick solution existence, uniqueness, and voltage balance guarantees without solving the full power flow problem. 


\section{Case Study}\label{sec:case}
This section first provides case studies to compare the approximations of the voltage balance requirement under the definition of $VUF_n$ and then demonstrates the results of the balancibility conditions using $PVUR$, $LVUR$, and $VUF_n$.   
 
 
\subsection{Tightness of VUF Approximations}
Here, we present instances of $C_{abc}$ and $r_{abc}$ to compare the results using approximation by bound \eqref{eq: vufrob_complex_conserv}, polytope approximation \eqref{eq: vufrobextre}, and LGR \eqref{eq: lgrfinal}. To get an estimate of the true optimal value of \eqref{eq: vufrob_bound}, we randomly generate $5\times 10^5$ points on $\partial\mathcal{V}^a_{in}\times\partial\mathcal{V}^b_{in}\times\partial\mathcal{V}^c_{in}$ and find the maximum value over the samples.\footnote{The sampling method is unsecure to directly use in the balancibility condition since it is a relaxation to \eqref{eq: vufrob_bound} and hence underestimates the voltage unbalance level.} For polytope approximation, we use $m=2,4,8,16,32$.  

\begin{figure}
\centering
\vspace{0cm}
\includegraphics[width=3.5in]{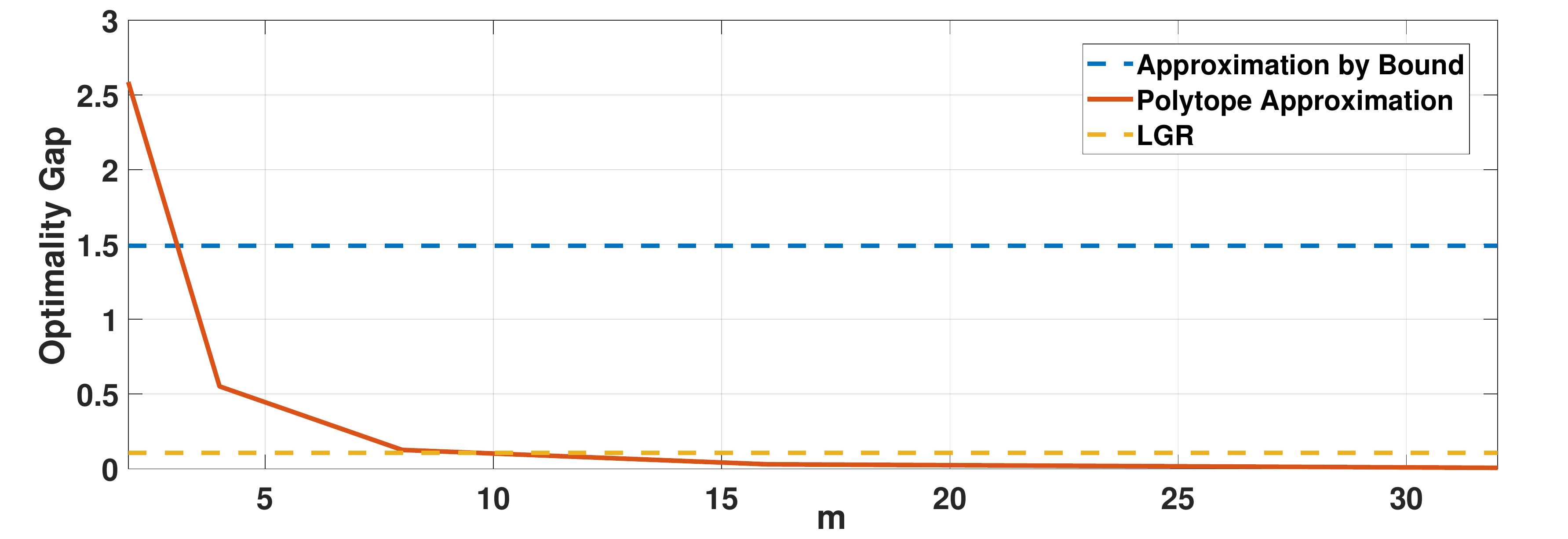}
\caption{Optimality Gaps of Approximation by Bound, Polytope Approximation, and LGR without Strong Duality}
\vspace{0cm}
\label{fig: fig1}
\end{figure}

As the first test, we select a center $C_a=(2,0)^{\top}$, $C_b=(-1,-\sqrt{3})^{\top}$, and $C_c=(-1,\sqrt{3})^{\top}$  (a balanced center for phases $a,b,c$); radius $r_a=r_b=r_c=0.6$; and tolerance $\epsilon=0.3$. Figure~\ref{fig: fig1} shows optimality gaps (i.e., the absolute difference from the estimate using sampling method) of the approximation by bound, LGR, and polytope approximations using different $m$. When $m$ is small, the polytope approximation has large optimality gap but the gap rapidly converges to zero as $m$ increases. LGR does not have strong duality (i.e., non-zero optimality gap) and hence gradually loses its advantage to the polytope approximation when $m$ increases. Approximation by bound provides an upper bound for the other approximation techniques. 

\begin{figure}
\centering
\vspace{0cm}
\includegraphics[width=3.5in]{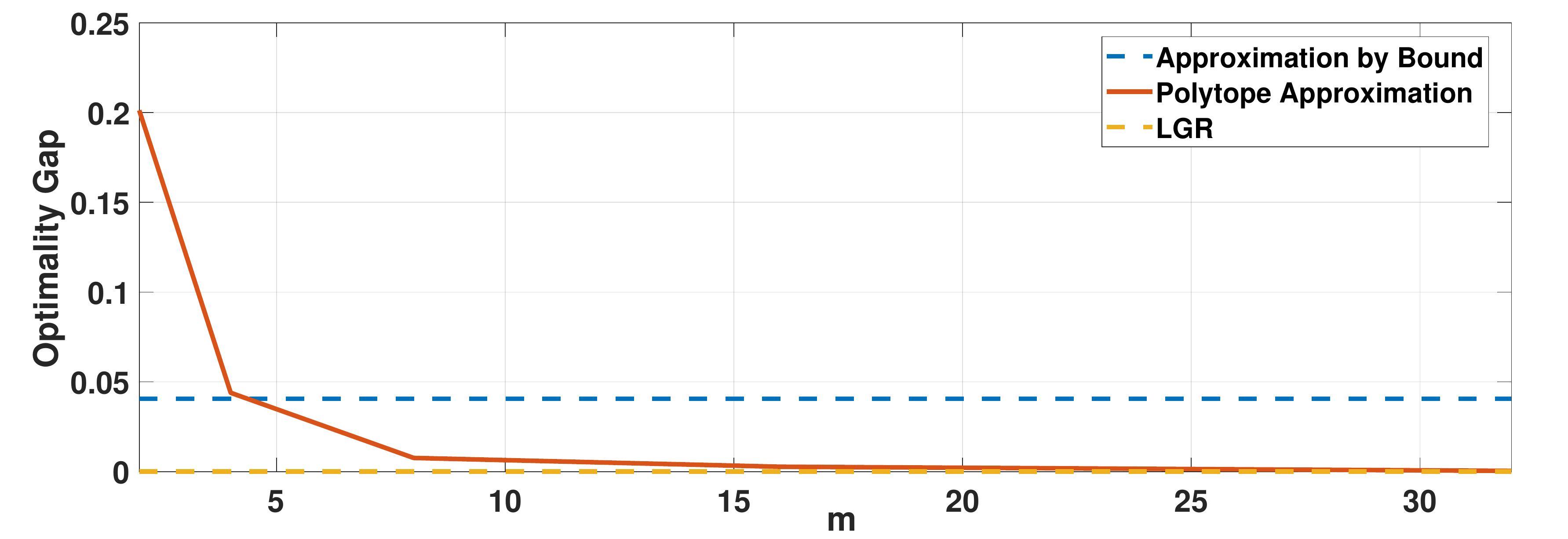}
\caption{Optimality Gaps of Approximation by Bound, Polytope Approximation, and LGR with Strong Duality}
\vspace{0cm}
\label{fig: fig2}
\end{figure}

As the second test, we select center $C_a=(3,0)^{\top}$, $C_b=(-1,-\sqrt{3})^{\top}$, and $C_c=(-1,\sqrt{3})^{\top}$; radius $r_a=r_b=r_c=0.1$; and tolerance $\epsilon=0.1$, which satisfy the sufficient conditions in Theorem~\ref{thm: suff1}. We then show the same groups of optimality gaps in Figure~\ref{fig: fig2}. Since strong duality holds between LGR and \eqref{eq: vufrob_bound}, we clearly see zero optimality gap in LGR. Meanwhile, the polytope approximation gradually converges to zero optimality gap as $m$ increases and the approximation by bound gives an upper bound. 

In practice, the tolerance $\epsilon$ is small $(\lesssim 5\%)$. Thus, when the set $\mathbf{U}_{in}$ (or, equivalently, $\mathcal{V}_{in}$) is small (i.e., the radius $r_{abc}$ is small), LGR has a higher chance of having a smaller optimality gap compared to the polytope approximation since Theorem~\ref{thm: suff1} is easier to satisfy.



\subsection{Balancibility Conditions}
To test the balancibility condition, we use a five-bus example system adopted from \cite{molly_5_bus} by only considering wye-connected PQ loads. As the nominal point $\mathbf{S}_L^0$, we choose a group of balanced loads on bus 4 ($10$ kW) and bus 5 ($50$ kW) and consider unbalanced loading at bus~$5$ via $\bm{\sigma}_L$ (i.e., incremental loads $\mathbf{S}_L-\mathbf{S}^0_L$). For $LVUR$ and $VUF_n$, we use the method of line-to-line voltage bounds \eqref{eq: llvbound} and LGR \eqref{eq: lgrfinal} respectively. As set up, we select $\bm{\sigma}_{L}^5=(10k,-5k,-5k)$ kW with $k=1,2,...,10$ and choose bus~$4$ as the critical node with respect to voltage balance. 

To compare the quality of the balancibility condition, we first seek the smallest tolerance $\epsilon$ such that the voltage balance requirements are satisfied over the $\mathcal{V}_L$ resulted from the solvability condition \eqref{eq:solv}. Then, we find the true unbalance by solving the power flow solutions with $\mathbf{S}_L$ and compute based on the unbalance definitions.

\begin{figure}[!t]
\centering
\vspace{0cm}
\includegraphics[width=3.5in]{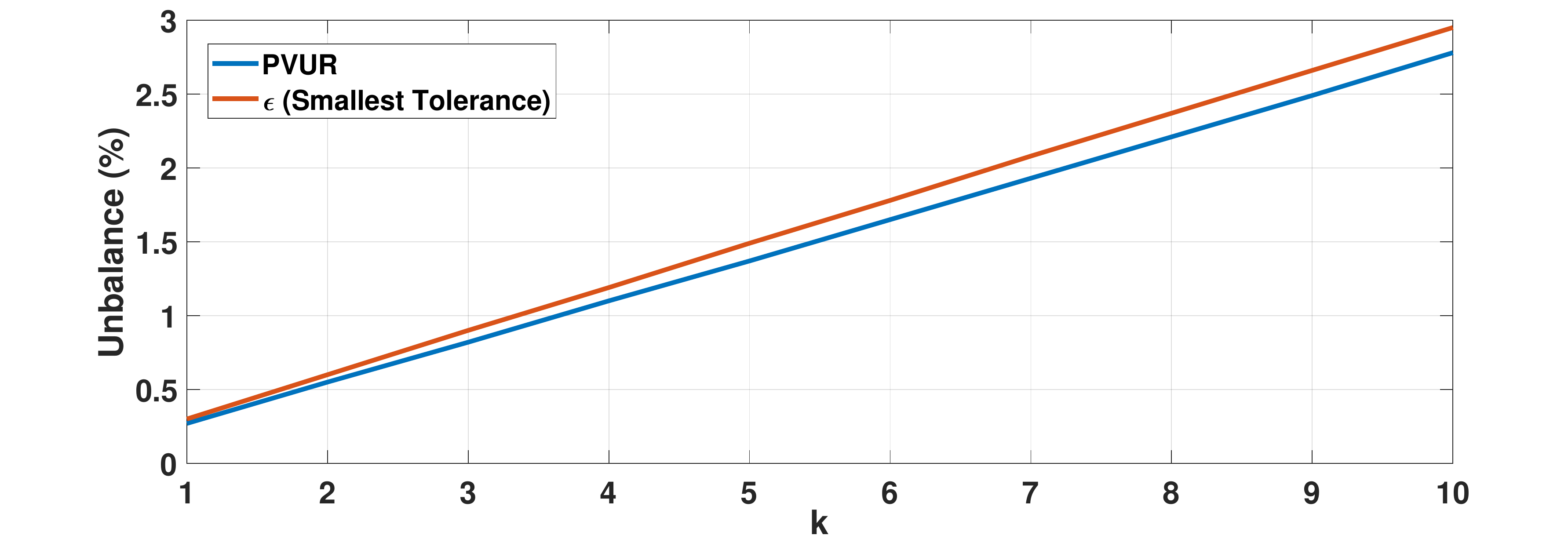}
\caption{Smallest Tolerance $\epsilon$ vs ``True'' PVUR under $\mathbf{S}_L$ from different~$k$}
\vspace{0cm}
\label{fig: fig3}
\end{figure}

\begin{figure}[!t]
\centering
\vspace{0cm}
\includegraphics[width=3.5in]{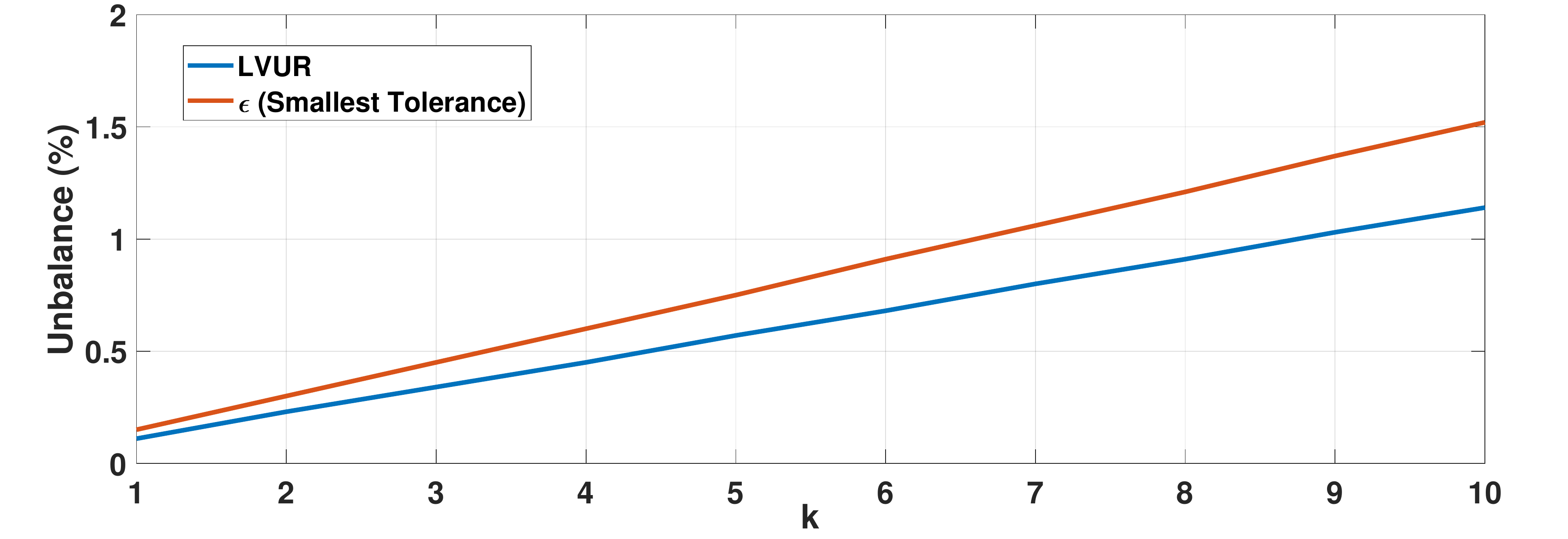}
\caption{Smallest Tolerance $\epsilon$ vs ``True'' LVUR under $\mathbf{S}_L$ from different~$k$}
\vspace{0cm}
\label{fig: fig4}
\end{figure}

\begin{figure}[!t]
\centering
\vspace{0cm}
\includegraphics[width=3.5in]{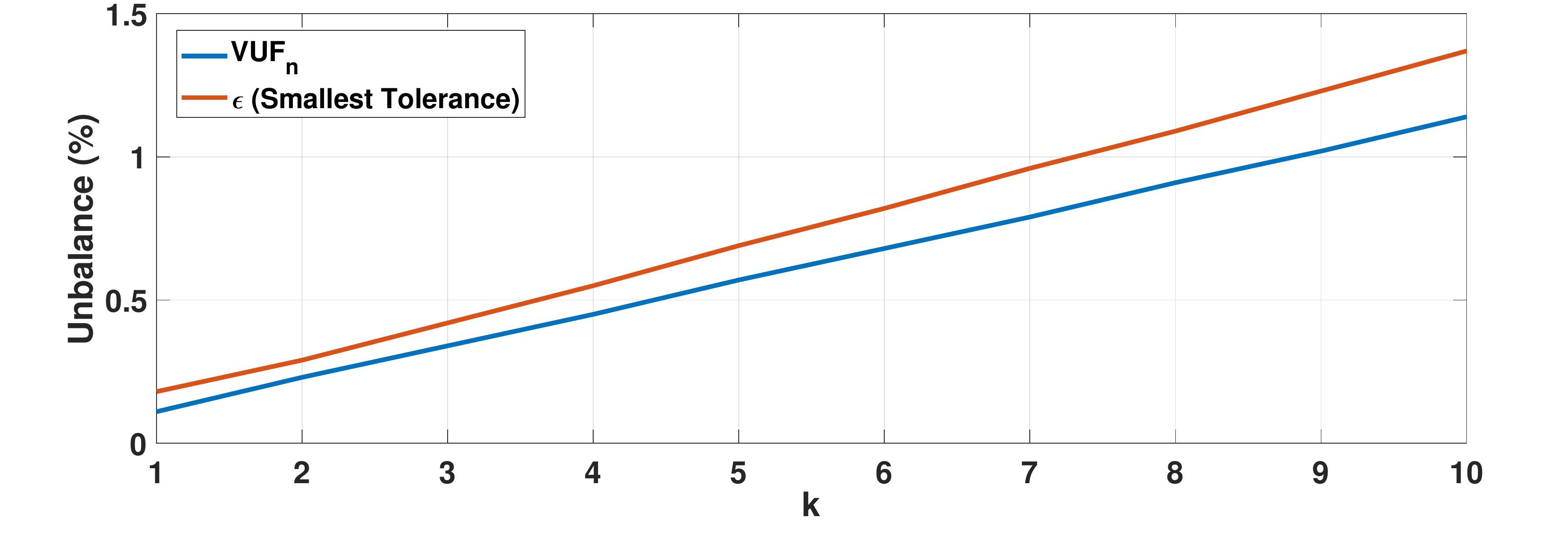}
\caption{Smallest Tolerance $\epsilon$ vs ``True'' VUF under $\mathbf{S}_L$ from different~$k$}
\vspace{0cm}
\label{fig: fig5}
\end{figure}

Figures \ref{fig: fig3}, \ref{fig: fig4}, and \ref{fig: fig5} show the true unbalance level and the smallest tolerance level $\epsilon$ under different unbalance definitions. As $k$ increases, the load on bus~5 becomes more unbalanced and hence induces a larger voltage unbalance level. We also see that there is a gap between $\epsilon$ and the true voltage unbalance level. The gap results from the fact that the balancibility condition provides voltage balance guarantees over $\mathcal{V}_L$, which contains the true solution $\mathbf{V}_L$ and hence a safe approximation. Since the gap is small, we conclude that the balancibility condition closely characterizes the system unbalance levels with the information of $\mathcal{V}_L$ without relying on the exact power flow solution $\mathbf{V}_L$. As~$k$ increases, $\mathcal{V}_L$ gets larger and hence results in a larger absolute gap. However, the ratio between $\epsilon$ and true unbalance level is approximately decreasing, which indicates reduced conservativeness. Further, we note that this conservativeness is system-dependent and can be improved if the solvability conditions (or general sets $\mathbf{U}_{in}$) are tighter.

\section{Conclusions and future work}\label{sec:conclu}
In this paper, we have proposed a concept called a \emph{balancibility condition} which combines the existing power flow solution tools (e.g., solvability conditions) with voltage balance requirements. The balancibility condition quantifies a power injection region which is guaranteed to have a unique and balanced power flow solution. We considered unbalance definitions from multiple organizations and derived closed-form reformulations or approximations to quantify the voltage unbalance level. We used a general model to describe the sets that contain the power flow solutions under uncertain power injections and gave theoretical guarantees on the quality of the approaches. As evaluation, we compared these approaches and demonstrated the benefits from the theoretical guarantees. We also compared the balancibility conditions associated with different unbalance definitions and demonstrated that the balancibility condition closely reflects the voltage unbalance level without excessive conservativeness.     

As future work, we will test all the balancibility conditions and approaches on more general and realistic distribution network models (e.g., delta-connected loads and ZIP loads). We will also improve the quality of the balancibility condition by optimizing the solvability condition together with the voltage balance requirements. Finally, we will consider relevant applications under uncertainty, including robust ACOPF problems.

\IEEEtriggeratref{12}
\bibliographystyle{IEEEtran}
\bibliography{ReferencesJournal1}

\end{document}